\documentclass[10pt,final]{siamltex}
\usepackage{amsmath}
\usepackage{amsfonts}
\usepackage{amssymb}
\usepackage{eucal}
\usepackage{cite}
\usepackage{graphicx}
\usepackage{color}
\usepackage[hypcap,font+=small,tableposition=top]{caption}
\usepackage[hypcap,labelfont=rm,font=normal,skip=0pt]{subcaption}
\usepackage[colorlinks,citecolor=blue!80!black,urlcolor=blue,linkcolor=red!80!black]{hyperref}
\usepackage{tikz}
\usepackage{pgfplots}
\usepgfplotslibrary{groupplots}
\pgfplotsset{compat=newest}
\usetikzlibrary{external}
\usepackage{cleveref}
\newlength\myindention
\DeclareCaptionFormat{parformat}%
{\hspace*{\myindention}#1#2#3}
\newcommand{\capsubref}[1]{\textup{(\subref{#1})}}
\newcommand{\subcap}[1]{\textup{(\subref*{#1})}}
\usepackage[colorinlistoftodos]{todonotes}
\makeatletter
\renewcommand{\listoftodos}[1][\@todonotes@todolistname]{\tikzexternaldisable\section*{#1}\hspace*{\fill}\newline\setcounter{tocdepth}{1}\@starttoc{tdo}\tikzexternalenable}
\renewcommand{\todo}[2][]{\tikzexternaldisable\@todo[#1]{#2}\tikzexternalenable}
\makeatother
\usepackage{forwarddefs}

\pgfplotscreateplotcyclelist{p4est2}{%
blue!80!black,solid,thick\\%
red!80!black,dashed,thick\\%
green!80!black,dotted,thick\\%
brown!80!black,solid,thick\\%
blue!80!black,dashed,thick\\%
red!80!black,dotted,thick\\%
green!80!black,solid,thick\\%
brown!80!black,dashed,thick\\%
blue!80!black,dotted,thick\\%
red!80!black,solid,thick\\%
green!80!black,dashed,thick\\%
brown!80!black,dotted,thick\\%
}

\pgfplotsset{every axis/.style={%
    title style       = {align = center, text width = 4.9cm},
    ymode             = log,
    xmin              = 0,
    xmax              = 60,
    xlabel            = Krylov iteration $j$,
    width             = 4.9cm,
    legend cell align = left,
    cycle list name=p4est2,
    scale only axis}}

\newlength{\convplotht}
\usepackage{multirow}

\def \addressices{Institute for Computational Engineering \& Sciences, The
  University of Texas at Austin, Austin, TX, USA}
\def \addressgeomech{Department of Geological Sciences and Department of
  Mechanical Engineering, The University of Texas at Austin, Austin, TX, USA}
\def \addresscims{Courant Institute of Mathematical Sciences, New York
  University, New York, NY, USA}

\begin{document}

\title{Solution of nonlinear Stokes equations \\ discretized by high-order
   finite
  elements \\ on nonconforming and anisotropic
  meshes, \\ with application to ice sheet
  dynamics\thanks{Support for this work was provided by
the U.S.\ Department of Energy Office of Science (DOE-SC), Advanced
Scientific Computing Research (ASCR), Scientific Discovery through
Advanced Computing (SciDAC) program, under award numbers
DE-FG02-09ER25914,
DE-11018096, and
DE-FC02-13ER2612,
and the U.S.\ National Science Foundation (NSF)
Cyber-enabled Discovery and Innovation (CDI) program under awards CMS-1028889
and OPP-0941678.  Allocations of computing time on TACC's Stampede under XSEDE, TG-DPP130002,
and on ORNL's Titan, which is supported by the Office of Science of the DOE under DE-AC05-00OR22725 are gratefully acknowledged.}}

\author{Tobin Isaac\footnotemark[2], Georg Stadler,\footnotemark[2] \footnotemark[3] \and
  Omar Ghattas\footnotemark[2] \footnotemark[4]}

\renewcommand{\thefootnote}{\fnsymbol{footnote}}
\footnotetext[2]{\addressices}
\footnotetext[3]{\addresscims}
\footnotetext[4]{\addressgeomech}
\renewcommand{\thefootnote}{\arabic{footnote}}

\newcommand{\firstblock}{(1,1)-block\xp}
\newcommand{\Firstblock}{(1,1)-Block\xp}
\newcommand{\secondblock}{(2,2)-block\xp}
\newcommand{\Secondblock}{(2,2)-Block\xp}

\maketitle

\begin{abstract}
  Motivated by the need for efficient and accurate simulation of the dynamics
  of the polar ice sheets, we design high-order finite element discretizations
  and scalable solvers for the solution of nonlinear incompressible Stokes
  equations. In particular, we focus on power-law, shear thinning rheologies
  commonly used in modeling ice dynamics and other geophysical flows.  We use
  nonconforming hexahedral meshes and the conforming inf-sup stable finite
  element velocity-pressure pairings $\tspace_\poly\times
  \tspace^\text{disc}_{\poly-2}$ or $\tspace_\poly \times
  \pspace^\text{disc}_{\poly-1}$, where $\poly\ge 2$ is the polynomial order
  of the velocity space.  To solve the nonlinear equations, we propose a
  Newton-Krylov method with a block upper triangular preconditioner for the
  linearized Stokes systems.  The diagonal blocks of this preconditioner are
  sparse approximations of the \firstblock and of its Schur complement. The
  \firstblock is approximated using linear finite elements based on the nodes
  of the high-order discretization, and the application of its inverse is
  approximated using algebraic multigrid with an incomplete factorization
  smoother.  This preconditioner is designed to be efficient on anisotropic
  meshes, which are necessary to match the high aspect ratio domains typical
  for ice sheets.
  As part of this work, we develop and make available extensions to two
  libraries---a hybrid meshing scheme for the \texttt{p4est} parallel adaptive
  mesh refinement library, and a modified smoothed aggregation scheme for
  PETSc---to improve their support for
  solving PDEs in high aspect ratio domains.  In a comprehensive numerical
  study, we find that our solver yields fast convergence that is independent
  of the element aspect ratio, the occurrence of nonconforming interfaces,
  and of the mesh refinement, and that depends only weakly on the polynomial
  finite element order.  We simulate the ice flow in a realistic description
  of the Antarctic ice sheet derived from field data, and study the parallel
  scalability of our solver for problems with up to 383 million unknowns.
\end{abstract}



\section{Introduction}
\label{sec:intro}

We design high-order finite element discretizations and scalable solvers for
incompressible nonlinear Stokes equations describing creeping flows of
power-law rheology fluids. Applications
include ice sheet dynamics
\cite{Hutter83}, mantle convection \cite{SchubertTurcotteOlsen01}, magma
dynamics \cite{McKenzie84} and other problems involving non-Newtonian
fluids \cite{GlowinskiXu11}.  Among the main challenges for the solution of
these problems are the presence of local features that emerge from the
nonlinear constitutive relation, the strongly varying and anisotropic
coefficients arising upon linearization, the incompressibility condition
leading to indefinite matrix problems, complex geometry and boundary
conditions, a wide range of length scales that may require highly-adapted
meshes with high aspect ratios, and large problem sizes that necessitate 
parallel solution on large supercomputers.
Our approach to cope with these challenges uses adaptively refined
nonconforming meshes, high-order inf-sup stable finite elements, and iterative
Newton-Krylov solvers combined with multilevel preconditioning techniques. We
focus in particular on the construction of efficient solvers and
preconditioners for discrete systems resulting from high-order
discretizations.

High-order finite element methods for partial differential equations (PDEs)
are attractive because, in many situations, the discrete solution converges
rapidly to the true solution as the approximation order $\poly$ is increased
or the characteristic mesh size $h$ is decreased.
However, this increased accuracy per degree of freedom compared to low-order
methods does not automatically translate into increased accuracy per unit of
computational work. This is due to the fact that matrices arising from
high-order discretizations are denser and, thus, more expensive to apply and
to solve systems with.
The cost of applying a matrix arising from a high-order discretization can be
reduced drastically if the work is shifted from memory operations to floating
point operations. This can be achieved using matrix-free implementations and
tensor-product approximation spaces and element operations on hexahedral
finite element meshes.  To precondition matrices arising from high-order
discretizations, low-order preconditioners based on the nodes of the
high-order discretization have proven efficient \cite{DevilleMund90,
  HeysManteuffelMcCormickEtAl05, Brown10, CanutoGervasioQuarteroni10,
  Olson07}.  These preconditioners allow fast construction and the use of
methods established for low-order discretizations.

Our approach to solving the nonlinear Stokes equations is an inexact
Newton-Krylov method, with a block preconditioning strategy for the linearized
equations, built from preconditioners for the \firstblock and for its Schur
complement. We consider a power-law rheology that involves
the second invariant of the strain rate tensor, for which the Newton
linearization results in a fourth-order aniso\-tropic tensor viscosity.
We pay particular attention to the interplay between discretization and
solver: seemingly minor differences in either the discretization or
the low-order preconditioner can vastly impact the performance of
conventional
solution methods for both diagonal blocks of the preconditioner.

Our driving application is the simulation of the dynamics of
continental-scale ice flows, which is a critical component of coupled climate modeling.
Predicting the contribution of ice sheets to sea-level rise is difficult
because of the complexity of accurately modeling ice sheet dynamics for the
full polar ice sheets and the large uncertainties in unobservable parameters
governing these dynamics \cite[Chapter 10, Appendix
6]{MeehlStockerCollinsEtAl07}.  To address these uncertainties, significant
effort has been focused on the development of inverse methods to infer ice
sheet model parameters from observations \cite{
PetraZhuStadlerEtAl12, MorlighemSeroussiLarourEtAl13}.  These inverse methods
require the repeated solution of ice flow equations for numerous parameter
fields that may vary over wide ranges, and many also require the repeated
solution of related adjoint ice flow equations. Hence, inverse methods
particularly stress the efficiency and robustness of solvers for nonlinear
Stokes equations.

A particular difficulty in ice sheet simulations is the high aspect ratio of
the computational domains, which is inherited by the discretization, leading
to anisotropic meshes. Discretizations with high-aspect ratio elements (and
problems with highly anisotropic material properties, which have many of the
same properties) are known to be challenging for implicit solvers and
preconditioners. The development of robust solvers for high aspect ratio
domains is also important in other Earth science and climate modeling
problems. In ocean flow models, for instance, three-dimensional implicit PDE
models are now being used \cite{KanarskaShchepetkinMcWilliams07}, whereas in
the past they were often replaced by two-dimensional approximations.




\paragraph{Related work}

Several recent articles develop scalable solvers
for Stokes problems with varying viscosity
\cite{BursteddeGhattasStadlerEtAl09, MayMoresi08, ElmanSilvesterWathen05,
GeenenurRehmanMacLachlanEtAl09, GrinevichOlshanskii09, FuruichiMayTackley11,
CaiNonakaBellEtAl14}.  These contributions use low-order stable or stabilized
finite elements for the discretization of the Stokes equations and address
nonlinearity mostly via
a Picard fixed point approach. Scalable solvers for high-order
discretizations of nonlinear scalar problems and extensions to linear 
incompressible flow problems are studied in
\cite{Brown10}. Various asymptotics-based approximations of the Stokes equations
are used for ice sheet and glacier modeling, which reduce the indefinite
Stokes equations to positive definite elliptic systems.  These simplifications
are justified by the large differences between horizontal and vertical
components of the velocity;
we refer to \cite{Hindmarsh04} for a comparison and discussion of the validity of
these different models. Ice sheet simulations using the full nonlinear Stokes
equations can be found, e.g., in \cite{GagliardiniZwinger08,
LengJuGunzburgerEtAl12, PattynPerichonAschwandenEtAl08}.

\paragraph{Contributions}

One of the main contributions of this paper is the design of discretizations,
solvers and preconditioners that allow the fast and scalable iterative
solution of nonlinear Stokes problems.  In particular, we obtain convergence
that, for a large class of realistic problems, is independent of the mesh
size, the presence of nonconforming interfaces in the mesh, and the element
aspect ratio, and depends only weakly on the polynomial order.  Another
contribution is the extension of low-order preconditioners for high-order
discretized problems to meshes with nonconforming interfaces and high aspect
ratio elements. 

In addition to analyzing our solver techniques on workstation-sized model
problems, we also demonstrate their performance and scalability on a series of
larger problems requiring a distributed memory parallel implementation,
including a simulation of the dynamics of the Antarctic ice sheet. The
simulation uses a geometry and temperature field derived from field data and
constitutes what we believe to be the first highly resolved nonlinear
Stokes-based continental scale simulations of the Antarctic ice sheet.

We have also developed publicly available tools for discretizing and solving
PDEs on high aspect ratio domains, such as (but not limited to) those occurring
for ice sheets.
One, presented in \cref{sec:space}, is an extension to the
\texttt{p4est} parallel adaptive mesh refinement (AMR) library \cite{BursteddeWilcoxGhattas11} that allows it
to construct favorable meshes for these domains. The other one, discussed in
\cref{sec:amg}, is a plugin for PETSc's
\cite{BalayBrownBuschelmanEtAl12}
generic algebraic multigrid (GAMG)
preconditioner to improve its effectiveness for these types of discretizations.



\paragraph{Limitations}
We use algebraic multigrid (AMG)
for preconditioning. An alternative would be to use geometric
multigrid (GMG), which builds a hierarchy through geometric coarsening of the
mesh. GMG can be tailored to only coarsen in certain directions while leaving
others unchanged (semicoarsening), which can be useful for anisotropic
geometries \cite{TezaurPeregoSalingerEtAl15}. We use the easier-to-use AMG,
in which we incorporate a minimal amount of geometric information.

Our simulations of ice sheet dynamics use a fixed temperature field and
geometry.  Simulations of evolving ice sheets would require coupling of the
nonlinear Stokes equations with a time-dependent advection-diffusion equation
for the evolving temperature, and with a kinematic equation for the evolution
of the ice sheet surface.  However, the solvers presented in this paper
carry over as important components in a time-stepping procedure for the
simulation of time-evolving nonlinear viscous flows.

\paragraph{Overview}%
This paper is organized as follows.  In \cref{sec:prob}, we discuss the form
of the nonlinear Stokes equations and boundary conditions that are the focus
of this work, their variational formulation, and their linearization.
In
\cref{sec:disc}, we present stable mixed-space finite element discretizations
and a discussion of adaptive mesh refinement.
In \cref{sec:linsol}, we give an overview of
our approach to solving the resulting discrete system of nonlinear equations
and the linearized counterparts.  The preconditioner for the \firstblock of
the linearized Stokes systems that arise at each Newton iteration, which is a
critical component of the linear solver, is presented in detail
in \cref{sec:fblock}, followed by a discussion on the preconditioner
for the Schur complement of the \firstblock in
\cref{sec:schur}.  We test our solver on a model geometry in
\cref{sec:nonlinex}, and then on a discretization of the Antarctic ice sheet
in \cref{sec:ant}, where we also study the scalability of our method.  We
conclude with a discussion in \cref{sec:discussion}.

\section{Nonlinear incompressible Stokes equations}
\label{sec:prob}

After specifying the Stokes equations with strain rate thinning
power-law rheology in \cref{sec:eqns}, we present the corresponding
variational form and argue existence of a unique solution in \cref{sec:weak}.
The linearization of the nonlinear equations is
presented in \cref{sec:nonlinsol}.

\subsection{Problem statement}
\label{sec:eqns}

On an open, bounded domain $\Dom\subset \mathbb{R}^3$ we consider the
incompressible Stokes equations
\begin{subequations}\label{eq:stokes}
\begin{align}
  -\Div \stress &= \vec{b}, &\xx&\in\Dom\label{eq:stokes1}, \\
  \div \uu &= 0, &\xx&\in\Dom\label{eq:stokes2},
\end{align} 
\end{subequations}
where $\uu$ is the flow velocity and $\vec{b}$ is a body force. The Cauchy
stress tensor $\stress$ depends on the strain rate tensor $\D(\uu) = \half
(\Grad \uu + {\Grad \uu}^\tr)$, its second invariant $\IIDu := \half
\frob{\D(\uu)}{\D(\uu)}$ and, possibly, other physical quantities such as a
temperature field. Here, ``$:$'' denotes the Frobenius product between
second-order tensors $\tens A = (\tens A_{ij})$ and $\tens B = (\tens B_{ij})$ defined by
$\tens A:\tens B = \sum_{i,j}\tens A_{ij}\tens B_{ij}$.
In the ice sheet problem, which is our driving
application, the stress tensor is given by Glen's flow law
\begin{equation}\label{eq:glen}
  \stress = -\pp \id + B(\T) (\IIDu + \varepsilon)^{\frac{1 - n}{2n}} \D\uu,
\end{equation}
where $\pp$ is the pressure, $B(\T)$ a positive-valued function of temperature
$\T$, $n\ge 1$ is the strain rate exponent, and $\varepsilon>0$ a small
regularization parameter that prevents infinite effective viscosity for $n>1$.
For $n=1$, \cref{eq:glen} reduces to a linear rheology, and it describes a
strain-rate weakening non-Newtonian fluid for $n>1$.  A common value used for
modeling the flow of glacial ice is $n=3$.  To complete the
definition of the boundary value problem \cref{eq:stokes}, it remains to
specify the boundary conditions.



In ice sheet simulations, different parts of the boundary require different
combinations of Dirichlet, Neumann, and Robin-type boundary conditions. This
makes ice flow a good problem for developing methods for other creeping flow
problems with complicated boundary conditions.  At the ice-air interface, the
homogeneous Neumann condition $\stress\nrml = \bs0$ holds. At the ice-water
interface, the normal stress matches the hydrostatic water pressure, i.e.,
$\stress\nrml\cdot\nrml = -\pp_{\text{w}}$ and the tangential components of
the boundary traction vanish.
At the base of the ice sheet, complex interactions occur between ice, water,
rock, and till.  In cold regions, the ice sticks to the bedrock, while in
temperate regions, water accumulates at the base and the
ice can slide but is subject to some amount of
friction.  A general way to describe these phenomena is to use a Dirichlet
condition in normal direction to describe melting and freezing at the base of
the ice sheet, combined with a Robin-type sliding law relating the tangential
component of velocity $\tang{\uu} = (\id -\nrml \otimes \nrml) \uu$ to the
tangential component of the stress through a function
$\beta(\cdot,\cdot,\ldots)$, i.e.,
\begin{equation}\label{eq:beta}
  \tang{\stress \nrml} = -\bta(|\tang\uu|,\T,\dots).
\end{equation}
Physically realistic descriptions of sliding must include the dependence of
\bta on $|\tang\uu|$ \cite{Zmitrowicz03}, which makes \cref{eq:beta} a
nonlinear boundary condition for the flow; in this work, however, we consider
the linear case $\bta(\tang\uu)=\bta(\xx)\tang\uu$, and no basal freezing or
melting. To summarize, we use the
following boundary conditions for the base \Bndr:
\begin{subequations}\label{eq:robin}
\begin{align}
  \tang{\stress \nrml} +\bta(\xx)\tang\uu  &= 0 &\xx&\in \Bndr,   \label{eq:robin1}\\
  \uu\cdot\nrml           & = 0  &\xx&\in \Bndr. \label{eq:robin2}
\end{align}
\end{subequations}


Polar ice sheets have a characteristic depth of less than 5 kilometers, while
they extend horizontally for thousands of kilometers.  Because of this
difference between length scales, modelers often simplify \cref{eq:stokes}
using asymptotic expansions that require assumptions about the
magnitude of the velocities and stresses in the ice sheet,
for instance, the shallow ice approximation \cite{Hutter83} and the
hydrostatic approximation \cite{Blatter95}.  The assumptions justifying these
simplifications do not hold for the entire ice sheet, which has led to
approaches that combine simplified models in the interior with Stokes
equations at outlet glaciers \cite{SeroussiDhiaMorlighemEtAl12}.  To avoid
these complications, we do not use simplified models and focus
on the efficient solution of the Stokes equations \eqref{eq:stokes} instead.

\subsection{Variational formulation}
\label{sec:weak}

Here, we define a variational form of \cref{eq:stokes,eq:robin} that defines
the fields $(\uu,\pp)$ as the unique solution in a 
vector
space $\cV \times \cM$ to
\begin{equation}\label{eq:weak}
  \int_{\Dom} \left[\visc(\uu)\D(\vv)\!:\!\D(\uu)
  - \pp \div\vv -\qq\div\uu\right]\ d\xx +
  \int_{\Bndr} \bta \tang\vv \cdot \tang\uu\ d\vec{s} = \ff(\vv),
\end{equation}
for all $(\vv,\qq)\in\cV \times \cM$, where $\visc(\uu) = \visc(\uu,\T) =
B(\T)(\IIDu + \varepsilon)^{\frac{1-n}{2n}}$, and \ff is the sum of the
effects of body and boundary forces.
We assume that $B(T)\in \cL^\infty(\Dom)$ is uniformly bounded from below,
that $\bta\in \cL^\infty(\Bndr)$ is nonnegative and that $\partial\Dom$ is
Lipschitz continuous. For simplicity, we assume that the Dirichlet boundary
conditions are homogeneous and are incorporated into the space $\cV$.

For a similar problem, Jouvet and Rappaz \cite{JouvetRappaz12} show that a
unique solution $(\uu,\pp)$ exists in the Dirichlet-conforming subspace of
$[\cW^{1,r}(\Dom)]^3\times \cL^{r'}(\Dom)$, where $r=1+1/n$ and $r'=1 + n$.
In \cref{sec:proof} we define a pair of spaces $\cV$ and $\cM$, in which
\cref{eq:weak} is well-posed. This pair is only slightly modified from the pair
above to account for the linear Robin boundary condition.

\subsection{Newton linearization}\label{sec:nonlinsol}

The Newton linearization of \cref{eq:weak} about a velocity-pressure pair
$(\uu,\pp)$ are equations whose solution $(\utl,\ptl)\in \cV \times \cM$
satisfies
\begin{equation}\label{eq:newton}
  \int_{\Dom} \left[\D(\vv)\!:\!(\dvisc(\uu)\D(\utl)) - \ptl \div\vv
    -\qq\div\utl\right]\
  d\xx + \int_{\Bndr} \bta \tang\vv \cdot \tang\utl\ d\vec{s} = -\vec{r}
  (\uu,\pp,\vv,\qq)
\end{equation}
for all $(\vv,\qq)\in \cV \times \cM$. Here, $\dvisc(\uu)$ is an
anisotropic 4th-order tensor given by
\begin{equation}\label{eq:rank4}
  \dvisc(\uu) = \visc(\uu) \id + \frac{\partial \visc(\uu)}{\partial \IIDu}
  \D(\uu) \otimes \D(\uu) =
  \visc(\uu)\left( \id  - \frac{n - 1}{2n}
    \frac{\D(\uu) \otimes \D(\uu)}{\IIDu + \varepsilon}\right),
\end{equation}
and $\vec{r}(\cdot)$ is the residual of \cref{eq:weak}. Here, ``$\otimes$''
denotes the outer product between two second-order tensors.
Compared to the Newton linearization~\cref{eq:newton}, the commonly used Picard
linearization of \cref{eq:weak}
neglects the anisotropic part of the fourth-order tensor $\dvisc(\uu)$.
Using a finite element discretization of \cref{eq:newton,eq:rank4} is only
marginally more complex than the Picard linearization, as the
action of $\dvisc(\uu)$ on $\D(\utl)$, which is all that is required in a
finite element method, can be computed using Frobenius products with trial and
test functions.  The operator $\dvisc(\uu)$ is also found in the adjoint
equations corresponding to \cref{eq:weak}, which are used in inverse methods
to infer uncertain parameters from observations \cite{PetraZhuStadlerEtAl12,
  MorlighemSeroussiLarourEtAl13}.

\section{Discretization}\label{sec:disc}

Our goal for discretizing \cref{eq:weak,eq:newton} is to obtain discrete
solutions that converge to the continuous solutions rapidly in terms of the
number of unknowns and time-to-solution.  Our approach combines locally
refined meshes with high-order finite element approximation spaces.  The
adaptive meshes we use are discussed in \cref{sec:space}, the finite element
approximation spaces for the Stokes equations are described in
\cref{sec:fedisc}, and computational aspects of the discretization are
addressed in \cref{sec:comput}.

\subsection{Meshing}\label{sec:space}

\begin{figure}\centering
  \tikzset{hybridpartition/.style={scale=4.0}}
  \input{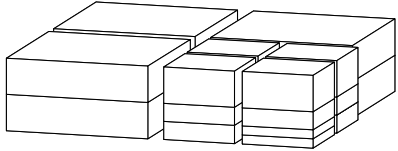}
  \caption{%
    Illustration of the type of mesh refinement provided by
    the \texttt{p6est} extension of the \texttt{p4est} parallel AMR library.
    The hexahedral elements are organized into columns, which are
    always assigned to a single process when the mesh is partitioned. Note that
    horizontal refinement is decoupled from vertical refinement.%
  }
  \label{fig:p6est}
\end{figure}

We use a hierarchical approach to mesh refinement, which starts with a coarse
mesh of conforming hexahedra.  This coarse mesh is expected to roughly
describe the geometry of the domain \Dom.  The fine mesh used for the finite
element discretization is obtained by hierarchical refinement of this coarse
mesh.  Refinement can be used to improve the resolution of the geometry of
\Dom or the resulting small hexahedra can simply inherit the geometry.  For
practical as well as numerical reasons, we require our refined meshes to obey
a 2:1 condition, where neighboring hexahedra can differ by only one level of
refinement, as illustrated in \cref{fig:hanging}.  We use the \texttt{p4est}
library to manage refinement, to enforce this 2:1 condition and to partition
the mesh between processes in parallel computations
\cite{BursteddeWilcoxGhattas11, IsaacBursteddeGhattas12}.  We use an extension
of the \texttt{p4est} library, developed by the first author, for meshing
three-dimensional problems in nearly two-dimensional domains.  This extension
uses a forest-of-quadtrees to manage independently-refinable columns of
hexahedra.  Meshes created with this extension have three key properties that
forest-of-octrees meshes lack.
\begin{enumerate}
  \item[(1)]
    Elements---and by extension degrees of freedom---are organized into
    columns that span the height of the mesh.  This organization can be
    exploited by solvers: see \cref{sec:amg}.
  \item[(2)]
    Each column of elements is assigned to a single process.  Split columns,
    which can appear in meshes created by the forest-of-octrees approach, can
    have negative consequences for solvers: see \cref{sec:amg}.
  \item[(3)]
    Hexahedra within a column may be locally refined in the vertical
    direction, allowing for more flexible refinement than the purely isotropic
    refinement of the forest-of-octrees approach.  We exploit this fine control
    over the element aspect ratio in our meshes of the Antarctic ice sheet:
    see \cref{sec:ant}.
\end{enumerate}
This type of mesh refinement is illustrated in \cref{fig:p6est}.  This
extension to the \texttt{p4est} library appears as collection of
``\texttt{p6est}''
data types and functions, so named because they build on elements of both the
``\texttt{p4est}'' interface for forests-of-quadtrees and the
``\texttt{p8est}'' interface for forests-of-octrees.  Documentation of the
\texttt{p6est} interface is available online.%
\footnote{For documentation of \texttt{p6est}, see
  \url{http://p4est.github.io/api/group__p6est.html}, as well as the
\texttt{example/p6est/test/test\_all.c} example distributed with the library,
which demonstrates the major I/O, mesh refinement, 
and visualization functions.}

Meshes for the simulation of ice flow must address the different length scales
inherent in the problem.  To accurately capture the vertical variations of the
state variables, a minimum vertical mesh resolution is necessary.  Most ice
sheet models use \abt 10 nodes in each vertical column and have a horizontal
resolution of 5 \km. Since the average thickness of the polar ice sheets is
\abt 2 \km, the width-to-height aspect ratio $\phi$ of these discretizations
is \abt 25.\footnotemark{} The Antarctic ice shelves, moreover, are typically
\abt 500~\m thick, so the same horizontal resolution results in $\phi \approx
100$.  We seek discretizations and methods for the Stokes equations that
support these aspect ratios.

\footnotetext{%
  In glaciology, one often uses the thickness-to-width aspect ratio
  $\epsilon=\phi^{-1}$ as the relevant
  limit in asymptotic expansions is $\epsilon\to 0$.  In this work, however,
  we prefer using $\phi$ because we consider the thickness of an ice sheet to
  be its characteristic length, in which we also measure its horizontal
  extent.%
 }



\subsection{Finite element discretization}\label{sec:fedisc}

In this section, we describe the finite element spaces used to discretize
the velocity and pressure spaces $\cV$ and $\cM$.

\subsubsection{Discrete velocity spaces}

\begin{figure}\centering
  \begin{minipage}[t]{0.49\textwidth}\centering
    \phantomsubcaption%
    \label{fig:hangingmesh}%
    \subcap{fig:hangingmesh}%
    \vtop{%
      \vspace{-0.6\baselineskip}%
      \hbox{%
        \includegraphics[height=0.45\textwidth]{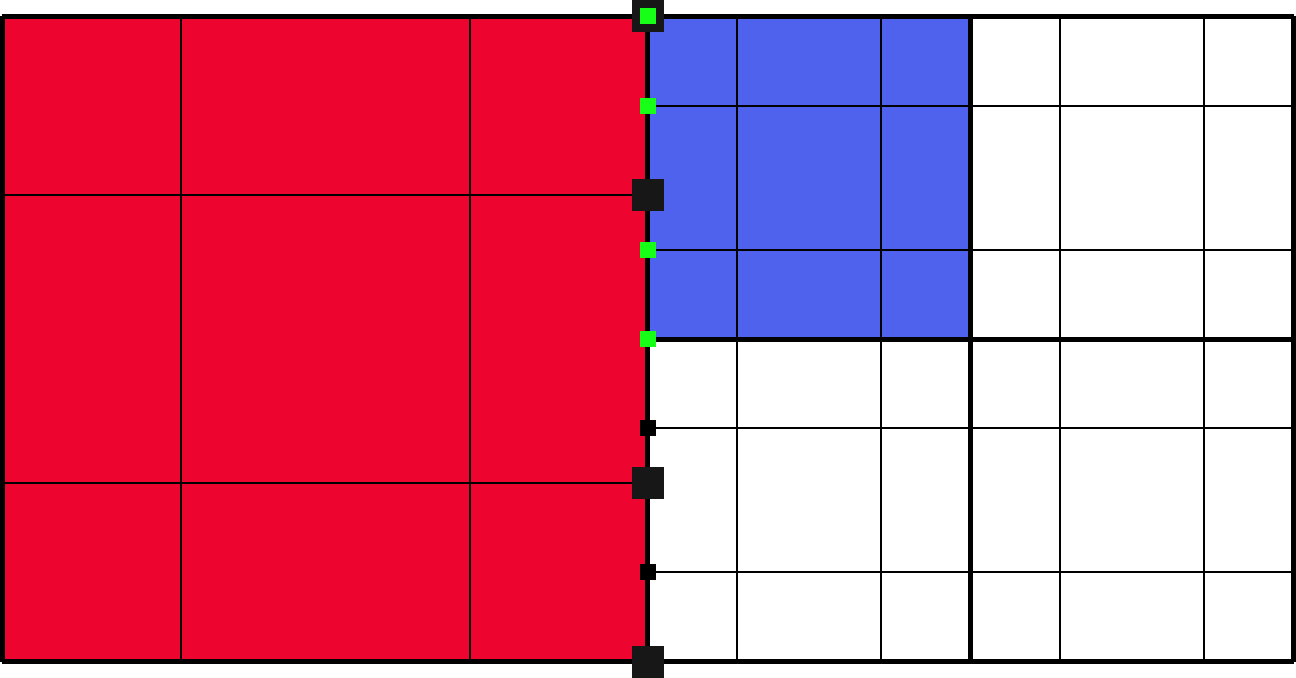}%
      }%
    }%
  \end{minipage}
  \begin{minipage}[t]{0.49\textwidth}\centering
    \phantomsubcaption%
    \label{fig:hangingnodes}%
    \makebox[0pt][l]{\subcap{fig:hangingnodes}}%
    \vtop{%
      \vspace{-0.6\baselineskip}%
      \hbox{%
        \tikzset{leftright/.style={scale=2.9},
          leftright/line/.style=ultra thick,
          leftright/bigline/.style=ultra thick}%
        \tikzset{
  leftright/.prefix style={},
  leftright/dot/.prefix style={draw=none,rectangle,inner sep=0pt,minimum size=4pt},
  leftright/bigdot/.prefix style={draw=none,rectangle,inner sep=0pt,minimum size=8pt},
  leftright/line/.prefix style={},
  leftright/bigline/.prefix style={thick},
}
\begin{tikzpicture}[leftright,
                    dot/.style=leftright/dot,
                    bigdot/.style=leftright/bigdot,
                    line/.style=leftright/line,
                    bigline/.style=leftright/bigline]

  \tikzstyle{connect}=[draw=gray,thin]
  \draw [color=red,thick,smooth,samples=200,domain=-1:1,bigline]
        plot (\x,{(\x+1)*(\x-1)*(\x-0.444721)});

  \draw [color=blue,smooth,samples=200,domain=0:1,line]
        plot (\x,{(\x+1)*(\x-1)*(\x-0.444721)});

  \draw (-1cm,0cm)           node[bigdot,fill=black] {}
        (-0.447cm,0.71554cm) node[bigdot,fill=black] {}
        (0.447cm,0.0cm)      node[bigdot,fill=black] {}
        (1cm,0cm)            node[bigdot,fill=black] {};


  \draw (0cm,0.44721cm)         node[dot,fill=green] {}
        (0.27639cm,0.15777cm)   node[dot,fill=green] {}
        (0.72361cm,-0.131672cm) node[dot,fill=green] {}
        (1cm,0.00cm)            node[dot,fill=green] {};


\end{tikzpicture}%
      }%
    }%
  \end{minipage}
  \caption{%
    \capsubref{fig:hangingmesh}
      A two-dimensional $\tspace_3$ mesh with a 2:1 nonconforming interface.
    \capsubref{fig:hangingnodes}
      The nodal values along a nonconforming interface.  Shown in green are
      the Gauss-Lobatto nodes of the smaller element, which do not align with
      those of the larger element: function values at these nodes depend
      on the values at the nodes of the larger element, so they are not
      included in the global nodal basis.  The matrix $R_i$ that
      interpolates a function to the nodes of the smaller element must
      interpolate cubic polynomials to the hanging nodes.  This
      polynomial interpolation is dense: the value at each of the hanging
      nodes is dependent on all of the independent nodes.  For
      two-dimensional nonconforming interfaces, $R_i$ is defined
      by tensor-product polynomial interpolation.%
  }%
  \label{fig:hanging}%
\end{figure}
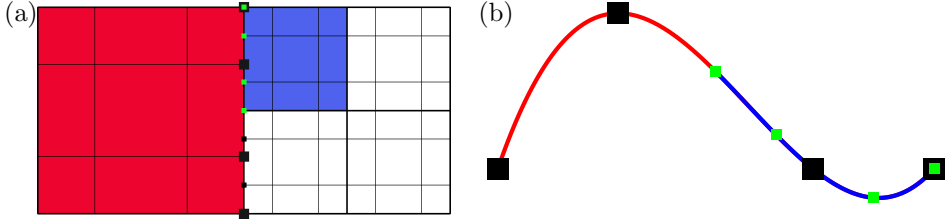

Given a mesh \cT of hexahedra $\{K_i\}$ with possibly nonconforming
interfaces, we define a finite-dimensional subspace of $\cC^0(\Dom)$ using
isoparametric Lagrange finite elements.  The $n_\poly := (\poly + 1)^3$ nodes
$\Xi_{\poly} = \{\bs\xi_{rst}\}_{0 \leq r,s,t \leq k}$
that define our Lagrange finite elements are the tensor-product Gauss-Lobatto
nodes of polynomial degree \poly on the reference domain $\refel = [-1,1]^3$.
These basis functions span $\tspace_\poly(\refel)$, the space of functions on
\refel that are the univariate polynomials of degree at most $\poly$ in each
of the coordinate directions.  We map \refel to an element $K_i$ by
$\varphi_i\in[\tspace_\poly(\refel)]^3$ and use the tensor-product Gauss nodes
of order $\poly$ for numerical quadrature.
We define the finite-dimensional space
\begin{equation*}
  V_{\cT,\poly} = \{v\in \cC^0(\Dom):\forall K_i\in\cT, v\circ \varphi_i \in
  \tspace_\poly (\refel)\},
\end{equation*}
and the velocity space $W_{\cT,\poly} = [V_{\cT,\poly}]^3 \cap \cV$.
For a conforming mesh, the set of element nodes $\Xi_\poly$ naturally defines
a set of global nodes $X_{\cT,\poly}$ for $V_{\cT,\poly}$ by the images of
element nodes. 
For a nonconforming mesh \cT, however, not all element nodes correspond to
global nodes, as shown in \cref{fig:hanging}.  To construct a set of global
nodes at a nonconforming interface, we thus ignore some element nodes: if an
element is more refined than its neighbor, the element nodes on the interface
with that neighbor do not contribute to the set of points that define the
global nodes.  Instead, these nodes are known as hanging or dependent nodes.

In general, function values at the nodes of element $K_i$ must be interpolated
from the global vector of nodal values by a restriction matrix $R_i$.  If
$K_i$ has no hanging nodes, then $R_i$ is simply a one-to-one association of
$\Xi_{\poly}$ to a subset of $X_{\cT,\poly}$; if $K_i$ has hanging nodes, then
$R_i$ interpolates values as described in \cref{fig:hanging}.  We use
identical trial and test spaces, so a global nodal matrix $A$ is assembled
from element nodal matrices $\{A_i\}$ by $A = \sum_i R_i^\tr A_i R_i$.

\subsubsection{Discrete pressure spaces}

We use inf-sup stable mixed finite element spaces to avoid the artificial
compressibility that can be introduced by stabilized discretizations of
incompressible flow. Additionally, to satisfy element-wise incompressibility,
we favor piecewise discontinuous pressure spaces $\cM_{\cT,\poly}$. This mass
conservation is
particularly important for ice sheet simulations, where the change of the mass
of the ice sheet is an important quantity of interest in climate projections.

The two most common choices for approximation on the reference cube are
$\pspace_{\poly-1}(\refel)$, which is the space of polynomials on \refel of
degree at most $(\poly-1)$, and $\tspace_{\poly-2}(\refel)$. We study two
possibilities for inf-sup stable velocity-pressure finite element spaces.
The pairing $\tspace_{\poly}(\refel) \times\pspace_{\poly-1}(\refel)$ has an
optimal order of convergence, and has an inf-sup constant
that is independent of $\poly$ and of the type of hierarchical local mesh
refinement we use \cite{HeuvelineSchieweck07}.  Its
inf-sup stability, however, degrades with increasing aspect-ratio $\phi$
\cite{AinsworthCoggins00}.  We find (see \cref{sec:schur}) that this
degradation can be significant for the element aspect ratios in our meshes.
An alternative pairing is $\tspace_{\poly} (\refel) \times \tspace_{\poly - 2}
(\refel)$, which has a suboptimal order of convergence, but its inf-sup
stability is uniform with respect to boundary layer refinement
making it appropriate for large values of
$\phi$ \cite{ToselliSchwab03}. For this pairing, the inf-sup constant decreases as $\cO(\poly^{-1})$;
however, for the moderate values of $\poly$ used in this work, this dependence
is not problematic.


\subsection{Computational aspects}\label{sec:comput}

For a mesh with $\nel$ elements, the number of degrees of freedom $\ndof$ in a
$\poly$-order finite element discretization is $\cO(\nel\poly^3)$, and the
number of nonzero entries in the matrix for a system of equations defined on
that space is $\cO(\nel\poly^6)=\cO(\ndof\poly^3)$.  This increasing density
means that sparse matrix-vector products based on globally assembled element
matrices are not efficient (in terms of memory operations) for large values of
$\poly$.  We therefore compute nonlinear residuals and apply linear operators
using a matrix-free approach to finite elements, where only the coefficients
and fields that define an operation are stored in memory and the operation's
application to a specific vector is assembled from all element contributions
at the time of application.  This approach requires $\mathcal{O}(\ndof)$
memory operations per matrix-vector product or residual calculation.  For
high-order elements, it is thus better suited to modern computer
architectures, where the bandwidth for memory operations is much narrower than
the bandwidth for floating-point operations.  This reduction in memory
operations comes at the expense of more floating point operations, but the
tensor structure of $\tspace_{\poly}(\refel)$ allows for all such element
computations to be reduced to repeated applications of one-dimensional compute
kernels, which can be highly optimized \cite{DevilleFischerMund02,Brown10}.
The same one-dimensional kernels are used to apply the restriction operators
$R_i$ for hanging nodes.


\section{Newton-Krylov method for nonlinear Stokes equations}\label{sec:linsol}

Our goal is to design a robust and scalable solver for the nonlinear Stokes
equations \cref{eq:stokes,eq:glen} with boundary conditions that
include \cref{eq:robin}. Ideally, the convergence should be
stable with respect to: (1) the element
size and the mesh refinement pattern, (2) parameters in the rheology
$\visc(\uu)$, (3) the Robin coefficient field \bta, (4) the polynomial order
\poly, and (5) the element aspect ratio $\phi$.  Here, we propose an inexact
Newton-Krylov solver, outlined in \cref{sec:disc_newton,sec:disc_krylov}
below. In \cref{sec:testprob}, we present a test problem that
allows us to study the effects of the factors listed above.  Variations of
this problem will be used in \cref{sec:fblock,sec:schur} to analyze and
optimize the convergence of our linear solver, and in \cref{sec:nonlinex}
for a nonlinear Stokes problem.

\subsection{Newton's method for nonlinear Stokes equations}\label{sec:disc_newton}

Given a velocity and pressure pair $(\uu,\pp)$, we (approximately) solve
\cref{eq:newton} for a search direction $(\utl, \ptl)$.  We then conduct a
line search in the direction $(\utl, \ptl)$ using the weak Wolf conditions
\cite{NocedalWright99} to guarantee that the nonlinear residual decreases.
Each Newton update is computed inexactly via a
Krylov-space iterative method, but with a tolerance that decreases in
subsequent steps so as to guarantee
quadratic convergence of the Newton
iterations \cite{EisenstatWalker96} close to the solution.  We ensure that discretization and
differentiation commute, so that the Jacobian obtained from discretizing
\cref{eq:newton} is equivalent to differentiating the discretization of
\cref{eq:weak}.
Discretization with one of the stable finite element pairs discussed in
\cref{sec:fedisc} results in a linear system with the symmetric saddle-point
system matrix
\begin{equation}
  A(\uu) = \begin{pmatrix} F(\uu) & B^\tr \\ B & 0 \end{pmatrix},
\end{equation}
where the \firstblock $F(\uu)$ is the discretization of the sum of the terms
involving \bta and $\dvisc(\uu)$ 
in \cref{eq:newton,eq:rank4}, and $B$ is the discretized divergence operator.


\subsection{Preconditioned Krylov method for linearized Stokes equations}
\label{sec:disc_krylov}

We solve systems involving $A(\uu)$ using preconditioned Krylov space methods,
typically restarted \text{GMRES}, or, if the preconditioner is not constant,
its flexible variant \text{FGMRES} \cite{Saad93}. As is well known, the
performance of Krylov methods critically depends on the availability of an
efficient preconditioner $\tilde{A}$ for $A(\uu)$. In the following, we use
the notation $A=A(\uu)$ and $F=F(\uu)$, i.e., in our notation we neglect the
dependence of $F$ and $A$ on $\uu$.
%
%
Due to the elliptic nature of $F$, a purely local preconditioner for $A$
cannot provide $h$-independent convergence and a multilevel preconditioner is
required. There are two main approaches for multilevel preconditioners for
incompressible flow problems, namely monolithic and block preconditioning
approaches.  The former approximates the saddle point system on each level of
a hierarchy and employs smoothers that are based on approximate local saddle
point solutions (i.e.,
Vanka-type smoothers) \cite{DamanikHronOuazzietAl09, Janka08}. This approach
typically requires a geometric mesh hierarchy or involves stabilized
discretizations.  In contrast, block preconditioners are built from
preconditioners for $F$ and for the Schur complement with respect to the
\firstblock, $S := -BF^{-1} B^\tr$. They allow to build on existing solvers
for elliptic systems and do not impose restrictions on the discretization
underlying $A$.
Due to this flexibility, we follow this latter approach and
use an upper-triangular block preconditioner $\tilde A$, such that the
preconditioned system becomes
\begin{equation}\label{eq:ps}
  A \tilde{A}^{-1} =
  \begin{pmatrix} F & B^\tr \\ B & 0 \end{pmatrix}
  \begin{pmatrix} \tilde{F} & B^\tr \\ 0 & \tilde{S} \end{pmatrix}^{-1} =
  \begin{pmatrix} I + \varepsilon_{\tilde{F}} & \varepsilon_{\tilde{F}}
  B^\tr
  \tilde{S}^{-1} \\
  B\tilde{F}^{-1} & I + \varepsilon_{\tilde{S}} \end{pmatrix},
\end{equation}
where $\varepsilon_{\tilde{F}} = I - F\tilde{F}^{-1}$ and
$\varepsilon_{\tilde{S}} = I -(-B\tilde{F}^{-1}B^\tr)\tilde{S}^{-1}$. Here,
the matrix $\tilde F$ is an approximation of the \firstblock $F$, and $\tilde
S$ is an approximation to the Schur complement $S$ with respect to the
\firstblock.
If $\varepsilon_{\tilde{F}}$ and $\varepsilon_{\tilde{S}}$ vanished,
$A\tilde{A}^{-1}$ would be a lower triangular matrix with all eigenvalues
clustered at 1, and a preconditioned Krylov method would converge in two
iterations. Thus, our target is to devise $\tilde F^{-1}$ and $\tilde S^{-1}$
such that $\varepsilon_{\tilde{F}}$ and $\varepsilon_{\tilde{S}}$ are small and $\tilde
A$-preconditioned Krylov methods converge quickly.  Additionally, it is
important that the setup time for the preconditioner $\tilde{A}$ is small,
because $\tilde{A}$ is recomputed as the Jacobian $F=F(\uu)$ changes.  In
\cref{sec:fblock,sec:schur}, we develop block preconditioners $\tilde{F}$ and
$\tilde{S}$ and study their properties.

In \cref{sec:nonlinex,sec:ant}, when we apply our preconditioner to problems
with dimensional coefficients, we combine the preconditioner $\tilde A$ that
we describe in this section with a diagonal rescaling of the linearized Stokes
system of the type described in \cite[Section 2.6]{MayMoresi08}.  This
rescaling tends to reduce the round-off errors that result from matrix
entries that differ by several orders of magnitude due to different units.

\subsection{Test problem setup}\label{sec:testprob}


Our test problems are based on problem C from the Ice Sheet Model
Intercomparison Project \cite{PattynPerichonAschwandenEtAl08}: \Dom is a
cutout of an ``infinite slab'', i.e., a sheet that is periodic in $x$- and
$y$-directions.  We use homogeneous Neumann boundary conditions on the top
surface of \Dom; on its base, we employ homogeneous 
Dirichlet conditions in the normal direction and Robin-type conditions in the tangential directions,
as in  \cref{eq:robin}.  We use the same boundary conditions when testing
systems in just the \firstblock $F$.  In linear test problems, we use
$\visc=1$ and $\bta=1$ for the material and
the sliding coefficients, unless specified otherwise.

\begin{figure}\centering
  \begin{minipage}[t]{0.43\textwidth}
    \phantomsubcaption
    \label{fig:stressxy}
    \makebox[\textwidth][c]{%
      \vtop{%
        \vspace{-0.6\baselineskip}%
        \hbox{%
          \includegraphics[width=\textwidth]{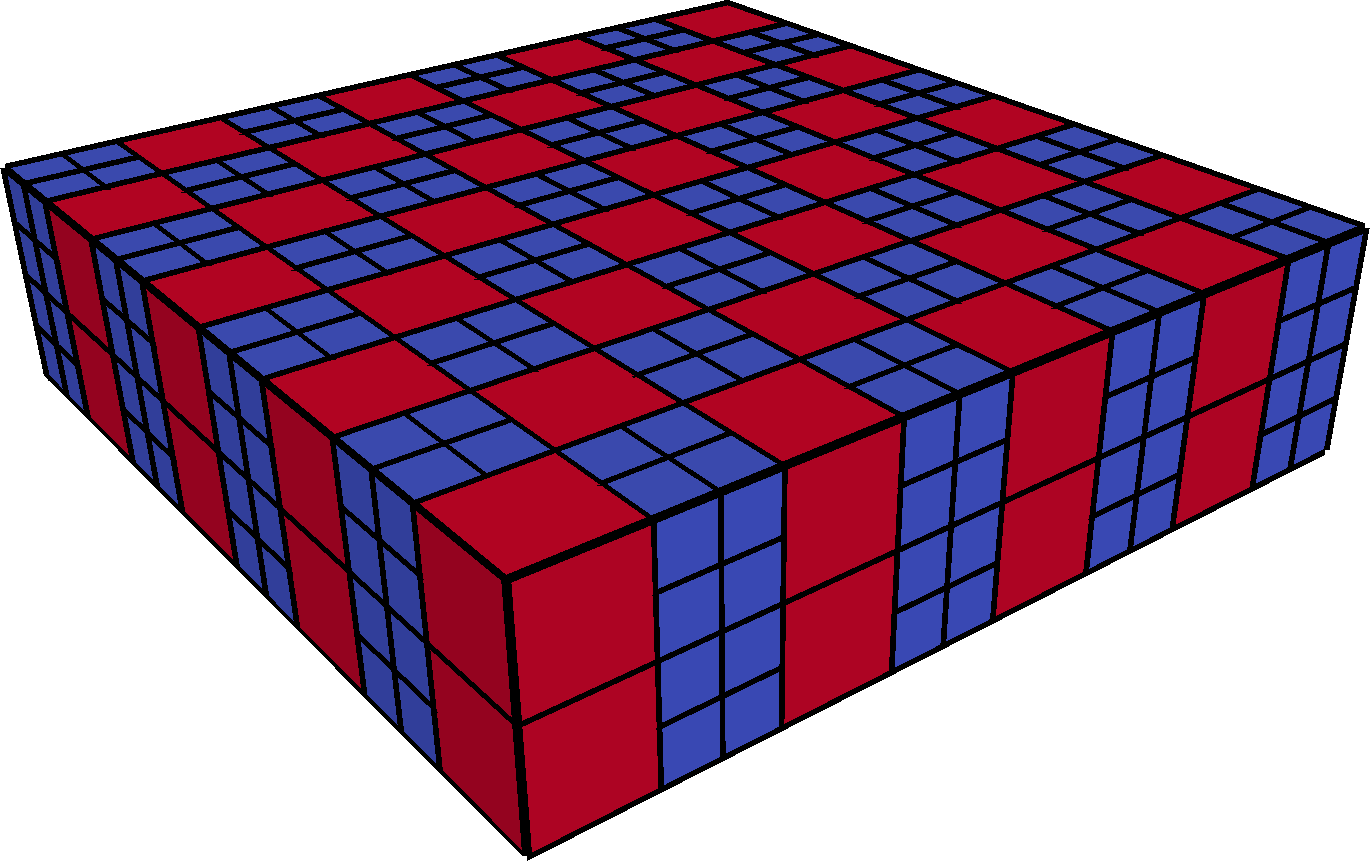}%
        }%
      }%
    }\hspace{-\textwidth}%
    \makebox[\textwidth][l]{\subcap{fig:stressxy}}%
  \end{minipage}
  \begin{minipage}[t]{0.43\textwidth}
    \phantomsubcaption
    \label{fig:stressmeshthin}
    \makebox[\textwidth][c]{%
      \vtop{%
        \vspace{-0.6\baselineskip}%
        \hbox{%
          \includegraphics[width=\textwidth]{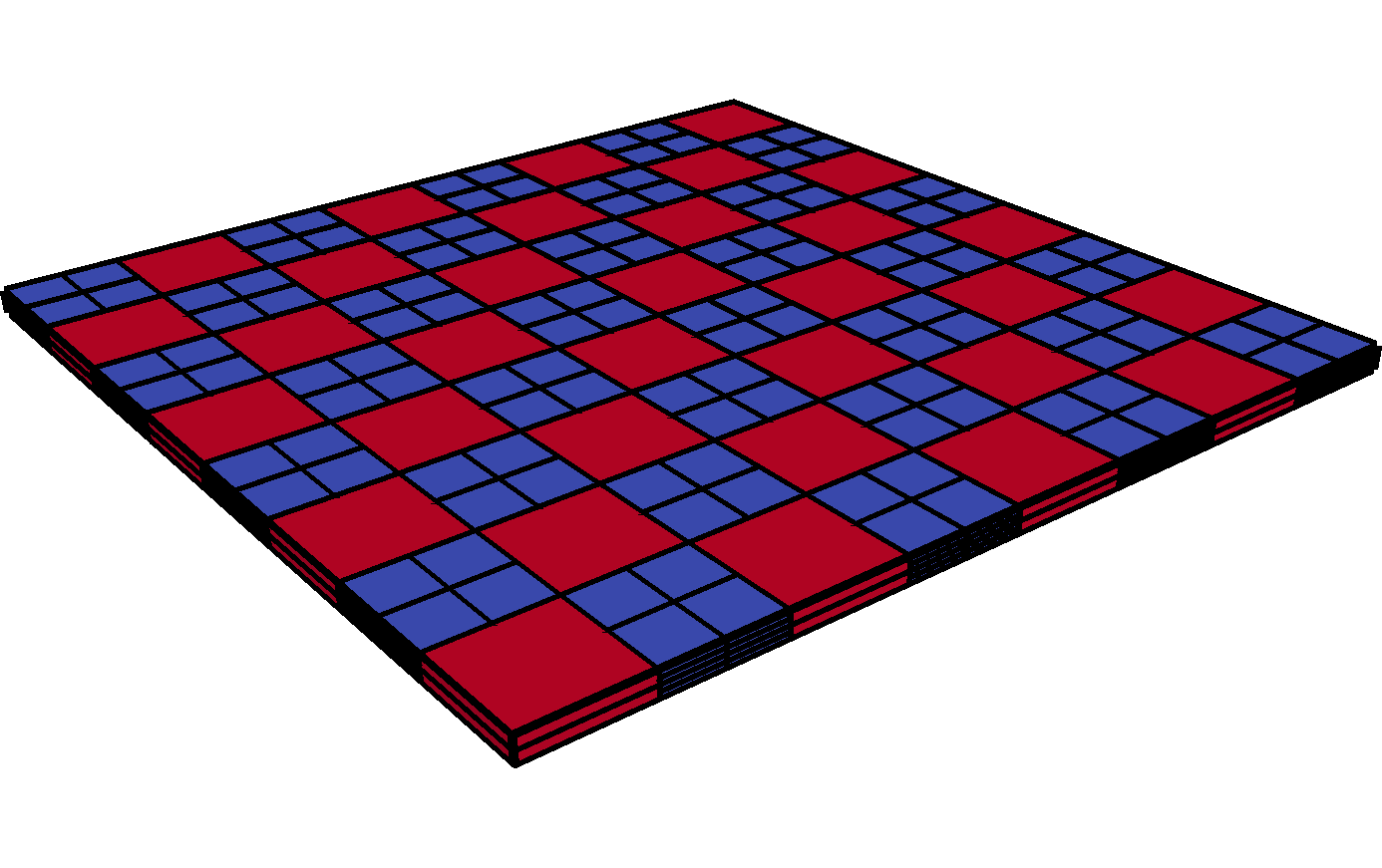}%
        }%
      }%
    }\hspace{-\textwidth}%
    \makebox[\textwidth][l]{\subcap{fig:stressmeshthin}}%
  \end{minipage}
  \begin{minipage}[t]{0.43\textwidth}
    \phantomsubcaption
    \label{fig:stressxyz}
    \makebox[\textwidth][c]{%
      \vtop{%
        \vspace{-0.6\baselineskip}%
        \hbox{%
          \includegraphics[width=\textwidth]{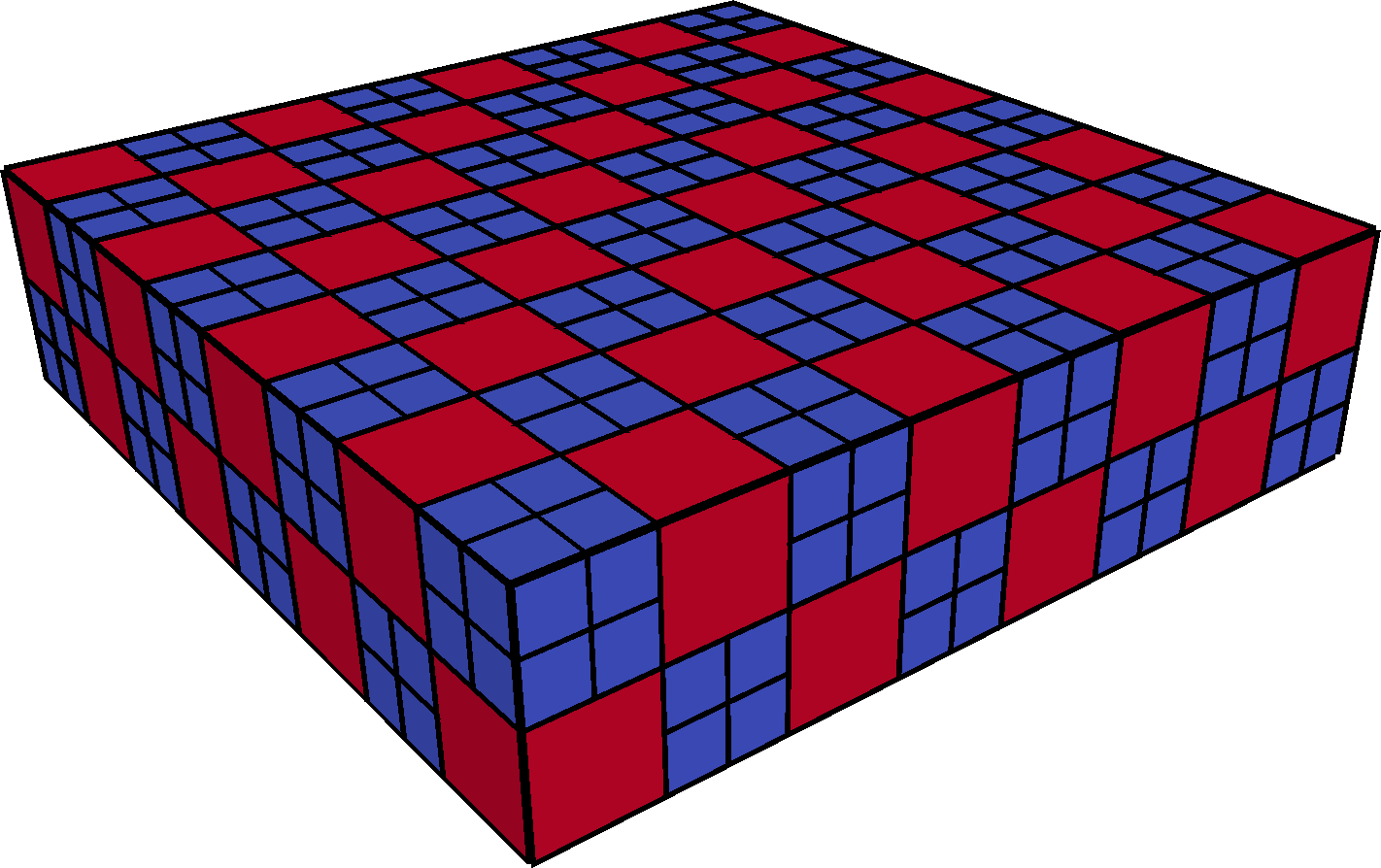}%
        }%
      }%
    }\hspace{-\textwidth}%
    \makebox[\textwidth][l]{\subcap{fig:stressxyz}}%
  \end{minipage}
  \begin{minipage}[t]{0.43\textwidth}
    \phantomsubcaption
    \label{fig:stresssol}
    \makebox[\textwidth][c]{%
      \vtop{%
        \vspace{-0.6\baselineskip}%
        \hbox{%
          \includegraphics[width=\textwidth]{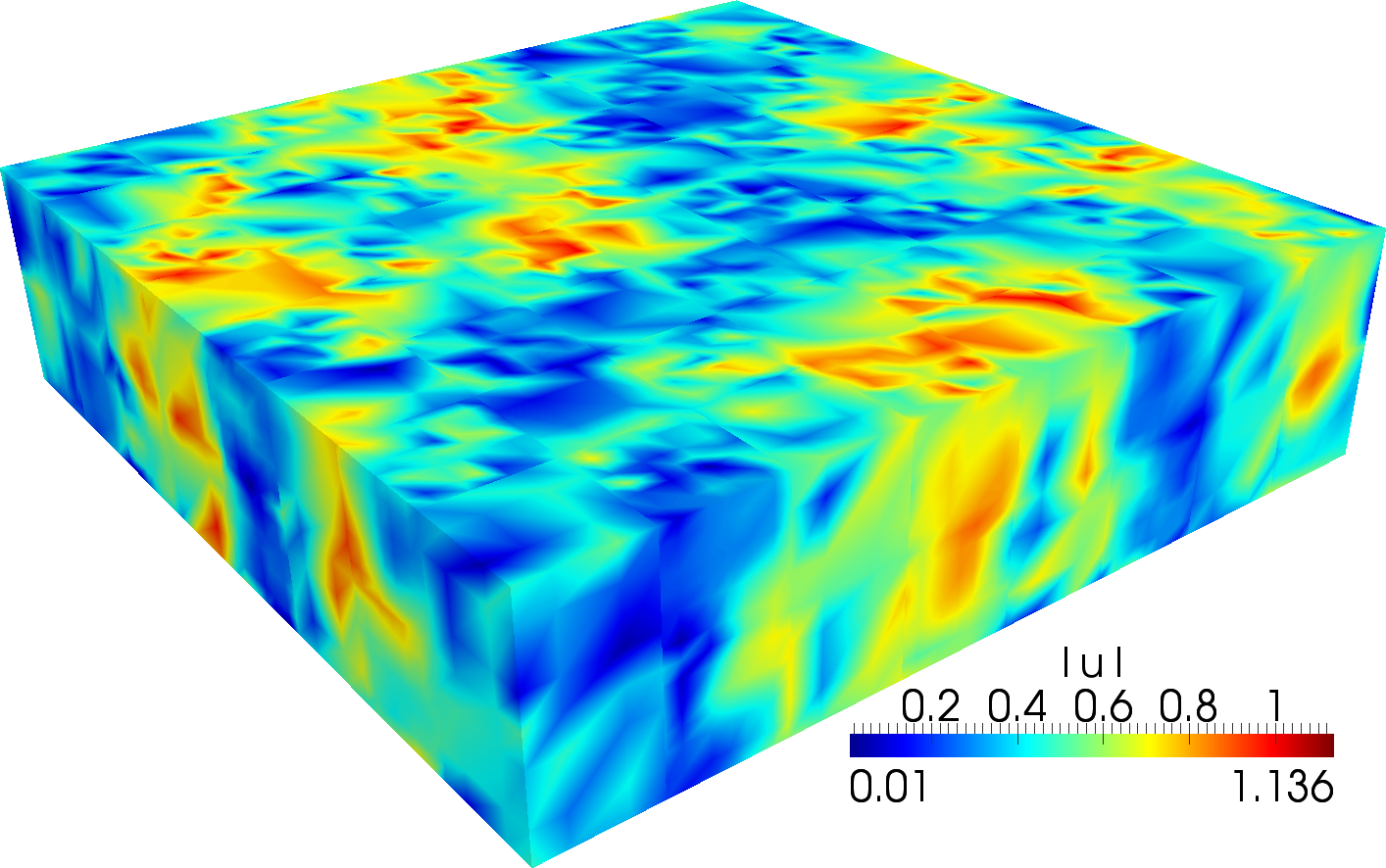}%
        }%
      }%
    }\hspace{-\textwidth}%
    \makebox[\textwidth][l]{\subcap{fig:stresssol}}%
  \end{minipage}
  \caption{%
    Meshes and manufactured solution used for the test problems.
    Elements are colored according to
    their level of refinement. 
      \capsubref{fig:stressxy}
        Mesh $\cT_{xy}$, which contains nonconforming faces in $x$- and
        $y$-directions, for $\phi=1$.
      \capsubref{fig:stressmeshthin}
        Mesh $\cT_{xy}$ for $\phi=10$.
      \capsubref{fig:stressxyz}
        Mesh $\cT_{xyz}$, which contains nonconforming faces in all directions.
      \capsubref{fig:stresssol}
        A manufactured solution $u^*$, constructed to contain variations with several
        length scales.%
  }%
  \label{fig:stressmesh}%
\end{figure}

To study the behavior of our solvers in the presence of nonconforming faces,
we use meshes for \Dom that have nonconforming faces throughout, as
illustrated in \cref{fig:stressmesh}.  We use two base meshes: $\cT_{xy}$,
which has nonconforming interfaces in only the $x$- and $y$-directions, and
$\cT_{xyz}$, which has nonconforming interfaces in all directions. Note that
the latter mesh cannot be decomposed into columns, so it represents a mode of
mesh refinement that appears when using octree-based meshes, but not when
using the hybrid AMR approach discussed in \cref{sec:space}. In all test
problems, \Dom has unit thickness, but we vary the length and width $L$ of the
domain.  The finite elements stretch with the domain, so that the number and
pattern of elements is the same, but the element aspect ratio $\phi$ varies
between 1 and 100; see \cref{fig:stressmeshthin}.  These meshes have 576
elements, and the number of degrees of freedom in the Stokes systems
discretized on these meshes is roughly $2300(k-1)k^2$.

For linear test problems, we test the effectiveness of our methods when
the residual contains multiple length scales.  To achieve this, we compute right hand
sides from  manufactured solutions, i.e., $b = A(\uu^*)(\uu^*,p^*)^\top$ when testing the
linear Stokes solver and $b = F(\uu^*)\uu^*$ when testing the \firstblock solver.
For that purpose, we create scalar-valued, spatially variable fields $s$ as
the sum of a Fourier series with random coefficients and pointwise random component:
%
\begin{equation}
  s(x,y,z) =\!\!\!\!\! \mathop{\sum_{j,k = 0}^N}_{(j,k)\neq(0,0)}\!\!\!\!
  (a^{j,k},b^{j,k},c^{j,k},d^{j,k})^\tr
  \begin{pmatrix}
    \cos\omega jx \cos \omega ky \\
    \cos\omega jx \sin \omega ky \\
    \sin\omega jx \cos \omega ky \\
    \sin\omega jx \sin \omega ky
  \end{pmatrix}
  |(j,k)|^{-\gamma} + e(x,y,z),
  \label{eq:randfield}
\end{equation}
where $N=10$, $\omega =2\pi / L$, $\gamma = 3/2$, and $i=1,2,3$. The
coefficients $a^{j,k}$, $b^{j,k}$, $c^{j,k}$ and $d^{j,k}$
are randomly chosen from  $[-1,1]$, but the Fourier coefficients
decay because of the $|(j,k)|^{-\gamma}:=(j^2+k^2)^{-\gamma/2}$ term, making $\gamma$ a control of
the smoothness of $s$.  The extra term $e$ is a random value from
$[-1/4,1/4]$, added at each node of the discrete vector.  To generate the
vector field $\uu^*$ that is the manufactured solution of our test problems for
the \firstblock solver, we generate a field $s$ for each component of $\uu^*$.
The magnitude of such a velocity field is shown in \cref{fig:stresssol}.  In
all our tests, the norm we use to report
the convergence of iterative solvers is the $\ell_2$-norm of the residual.

\section{Preconditioning the \firstblock}\label{sec:fblock}

The \firstblock $F$ occurring in the Stokes system is similar to the operator
arising in linear elasticity. If we neglect boundary conditions, its nullspace
is given by the rigid-body modes.
To approximately apply the inverse of $F$, we use algebraic multigrid (AMG)
and, in particular, we use the smoothed aggregation multigrid (SA-AMG), which
has theoretically proven convergence bounds \cite{VanekBrezinaMandel01}.
SA-AMG uses Galerkin projections to create coarse approximation spaces, where
the coarse space is embedded in the fine space by a prolongation matrix $P$.
$P$ is constructed by first creating a projector $\tilde{P}$ that projects
coarse ``aggregate'' nodes onto disjoint sets of fine nodes, followed by
creating $P$ from $\tilde{P}$ by applying a local smoothing operation based on
$F$, while ensuring that the near nullspace of $F$ (the nullspace in the
absence of boundary conditions) is well-approximated in the coarse space
\cite{VanekMandelBrezina96}.
To construct the prolongation and coarse matrices of the hierarchy, we use
SA-AMG as implemented by PETSc's GAMG preconditioner.  To build smoothers for each
level, we use PETSc's {\tt {KSP}} and {\tt PC} objects for defining Krylov
methods and preconditioners.
In this framework, we have three main design parameters, namely:
\begin{enumerate}
    \item[(1)]
      The matrix $\tilde{F}$ used to construct the multigrid hierarchy. As
      multigrid is only used as preconditioner and not as solver, $\tilde F$
      can be based on a lower-order (and thus sparser) approximations of the
      high-order element discretization of the \firstblock used in the
      residual computation.
    \item[(2)]
      How the hierarchy is coarsened.  If fine nodes are only grouped into
      aggregates when they are strongly connected to each other, then there
      will be more aggregates and a less aggressively coarsened hierarchy.  In
      standard SA-AMG, connections between fine nodes are ignored when forming
      aggregates if they are weaker than some threshold $\theta$.
    \item[(3)]
      The choice of smoothing operation on each level of the hierarchy.
\end{enumerate}
These three choices cannot be made independently.  The aggregation strategy
and the smoothers are related because the error components not damped by the
smoother must be corrected on the coarser levels.  Moreover, the choice of
$\tilde{F}$ affects the hierarchy and the effectiveness of different
smoothers.  In this section, we study these parameters on linear systems
involving the \firstblock $F$ rather than the full Stokes system $A$.

We note that two additional parameters, the size of the coarsest grid and the
solver used for the corresponding system, can also impact the scalability and
performance of multigrid preconditioners.  In this work, however, we have used
the default coarse grid size of GAMG (which is less than 50 unknowns), and
we use a direct solver for the corresponding system in all of our numerical
experiments. Note that in parallel simulations, we repartition the coarsened
multigrid levels to subsets of processors, and, in particular, the system for
the coarsest grid only uses a single processor.
%


\subsection{SA-AMG aggregation and smoothers}\label{sec:amg}

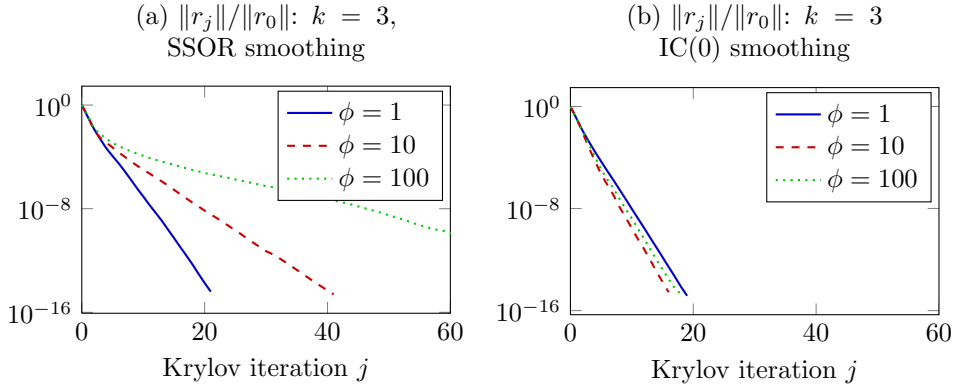
\begin{figure}\centering
  \begin{minipage}[t]{0.49\textwidth}\centering
    \phantomsubcaption%
    \label{fig:sor}%
    \begin{tikzpicture}[baseline]
      \begin{axis}
        [
          title = {\subcap{fig:sor} $\|r_j\|/\|r_0\|$: $\poly=3$,\\ SSOR smoothing},
        ]

        \foreach \i in {1,10,100} {%
          \addplot table
            [
              x index = 0, 
              y index = 11,
            ]
            {figures/sor/sor.\i.monitor};
            \addlegendentryexpanded{$\phi=\i$}
        }
      \end{axis}
    \end{tikzpicture}
  \end{minipage}
  \begin{minipage}[t]{0.49\textwidth}\centering
    \phantomsubcaption%
    \label{fig:stressbase}%
    \begin{tikzpicture}[baseline]
      \begin{axis}
        [
          title = {\subcap{fig:stressbase} $\|r_j\|/\|r_0\|$: $\poly = 3$\\ IC(0) smoothing},
        ]

        \foreach \i in {1,10,100} {%
          \addplot table
            [
              x index = 0, 
              y index = 11,
            ]
            {figures/base/base.\i.monitor};
            \addlegendentryexpanded{$\phi=\i$}
        }
      \end{axis}
    \end{tikzpicture}
  \end{minipage}
%
  \caption{%
    Convergence histories for the \firstblock test problem described in
    \cref{sec:testprob}, discretized with $\poly = 3$ on mesh $\cT_{xy}$
    (\cref{fig:stressxy}) at different aspect ratios.  The solver is GMRES
    preconditioned by one SA-AMG V-cycle. The solver in \subcap{fig:sor} uses
    a Chebyshev(1)/block-Jacobi/SSOR smoother, i.e., one optimally-damped
    application of block-Jacobi/SSOR; the solver in
    \subcap{fig:stressbase} is the same, except that it uses an incomplete
    Cholesky factorization (IC(0)) instead of SSOR.%
  }%
  \label{fig:base}
\end{figure}

The convergence of multigrid using pointwise smoothers such as Jacobi and
Gauss-Seidel is known to degrade for anisotropic problems when isotropic
coarsening is used \cite[Chapter~4]{TrottenbergOosterleeSchuller01}.  This
slowdown can also be seen in \cref{fig:sor}, where, for different element
aspect ratios $\phi$, we show the convergence based on a symmetric
Gauss-Seidel (e.g., symmetric successive over-relaxation, SSOR) smoother for
hierarchies created with aggregation threshold $\theta = 0$. This behavior
occurs because for large $\phi$, horizontal variations in the error have less
energy in the operator norm than vertical variations, and are thus not reduced
sufficiently.  There are two primary approaches to addressing this issue.  One
can use pointwise smoothers and use smaller aggregates to construct an
operator hierarchy such that the undamped error components are well
represented on the coarser meshes and can be corrected there (semicoarsening),
or one can continue to use full-sized aggregates and use a smoother that more
effectively dampens all error components (non-pointwise smoothing).  

Although some form of semicoarsening could be effective for our problems, we
have not had success with it. In standard SA-AMG, semicoarsening is
accomplished by only aggregating strongly connected nodes, for which the
corresponding matrix entries satisfy 
\begin{equation}\label{eq:threshold}
|F_{ij}|^2\geq \theta |F_{ii}||F_{jj}|.
\end{equation}
While this heuristic is generally good for scalar elliptic problems, for the
vector system $F$ it is problematic.  It is known that for a fixed threshold
$\theta > 0$, a sufficiently large aspect ratio $\phi$ will cause the standard
SA-AMG algorithm to break down for operators whose nullspace contains the
rigid-body modes \cite[Section 5.3]{GeeHuTuminaro09}.  Even before this
breakdown, we observed poor convergence rates with $\theta > 0$ and pointwise
smoothers.  A better heuristic for anisotropic problems than
\cref{eq:threshold}, which is reported to result in good convergence for
anisotropic elasticity problems, can be found in  \cite[Section
4]{KarerKraus10}, but even in this case, semicoarsening produces larger coarse
operators and deeper hierarchies, which increase the cost of storing and
applying the preconditioner.  We therefore choose not to use thresholding when
constructing our hierarchies, i.e., we do not consider the magnitudes of
matrix entries when constructing our aggregates, only the nonzero pattern of
the matrix.  This choice results in aggressively coarsened hierarchies with
smaller operator complexities, and has the additional advantage that it allows
us to reuse the projection matrices for multiple Newton steps, reducing the
cost of the subsequent preconditioner setups.

Since we must select a smoother that is compatible with our aggressive coarsening,
we use approaches based on incomplete factorizations of the matrix
$\tilde{F}$. In parallel, we only compute the incomplete factorization for
each process's diagonal block of the distributed $\tilde{F}$ matrix, which
amounts to a block-Jacobi smoother.  To build a stationary smoother from this
block-Jacobi/IC(0) smoother, an estimate of the largest eigenvalue of the
smoothed operator is required, either to calculate a single damping parameter
or the coefficients of a Chebyshev polynomial.  We numerically estimate the
largest eigenvalue using an iterative method, which adds to the setup cost of
the multigrid preconditioner.  To avoid this setup cost, one can alternatively
use a non-stationary smoother that does not require damping, such as a few
iterations of a Krylov method. This is a good choice when only a small
reduction in the residual is required, as for instance in the early iterations
of an inexact Newton-Krylov method.

We note that in certain parameter regimes---large values of $\beta$ or no-slip
Dirichlet conditions, and a fixed horizontal resolution---incomplete
block-Jacobi factorizations have been observed to work well as stand-alone
preconditioners, even with fine vertical resolution \cite[Section
3.1]{BrownSmithAhmadia13}. However, we are not aware of a systematic
study of the parameter regimes where this good performance can
be observed.

Using local Fourier analysis, one can show that geometric multigrid smoothed
by incomplete factorization can give $\phi$-independent convergence for scalar
elliptic problems \cite[Chapter 7]{TrottenbergOosterleeSchuller01}.  A key
component of this analysis is the ordering of degrees of freedom in the
incomplete factorization: tightly-coupled degrees of freedom within a column
must be sequential.  We can order the degrees of freedom in our original
system this way, but the standard techniques for generating aggregates, such
as GAMG's default method of randomized maximal independent set (MIS) selection
\cite{Adams98}, result in coarse discretizations having no special spatial
structure. 

Because of this, the multigrid hierarchies created by randomized MIS
aggregation are efficient for solving the \firstblock equations at a high
aspect ratio only if the boundary conditions are strong enough, which we
demonstrate in \cref{fig:basehagg}.  If the basal boundary has full Dirichlet
conditions, or a strong Robin coefficient $\beta$ ($\beta \geq 1$ in our
non-dimensional model problem), we see that GMRES preconditioned by such a
V-cycle converges with near $h$-independence.  If, however, we have
$\beta\approx 0$, then the V-cycles become less efficient and lose
$h$-independence.  The culprits for this loss of efficiency are error modes
which are highly oscillatory in the $x$- and $y$-directions but nearly constant
in the $z$-direction.  At high aspect ratios, the stress energy of these modes
is low in the absence of a strong boundary condition, so the smoother does
little to dampen them.  If the coarse mesh correction is poor, then it will
introduce error components in these modes.

Our solution to this problem of weak boundary conditions is to use meshes
whose degrees of freedom form columns such as the hybrid AMR meshes discussed
in \cref{sec:space}, and to
modify the aggregation technique used by GAMG in order to preserve the
column-structure of these degrees of freedom. We first construct
aggregates of columns using a standard aggregation technique such as
randomized MIS, 
and then partition each aggregated column into the final node
aggregates.  When subdividing a column, we try to ensure that each node
aggregate is at least three nodes tall, because three is usually the minimum
diameter of aggregates used in standard SA-AMG.  The multigrid hierarchies
created by this method are much closer to achieving $h$-independent
convergence for weak boundary conditions, as we demonstrate in
\cref{fig:basehdofcol}.  We note that this method appears to be effective even
when the number of degrees of freedom in each column varies as they do in our
test meshes. 
We also note that this technique does not depend on our particular choice of
hexahedral elements, but could also be applied, e.g., to meshes with
triangular prism elements, which are commonly used in modeling ice sheets and
other high aspect ratio domains.  We implement this aggregation technique in
the DofColumns plugin\footnote{The DofColumns plugin is available from
  \url{https://bitbucket.org/tisaac/dofcolumns/}.} for GAMG.

\begin{figure}\centering
  \renewcommand{\arraystretch}{1.1}
  \begin{minipage}[t]{0.49\textwidth}\centering
    \phantomsubcaption%
    \label{fig:basehagg}%
    \begin{tabular}{|c||c|c|c|c|}
      \multicolumn{5}{c}{\subcap{fig:basehagg} standard MIS} \\ \hline
      $\ell$ & $\beta=1^0$ & $1^{-2}$ & $1^{-4}$ & $1^{-8}$ \\ \hline
      0      & $0.14$      & $0.14$   & $0.57$   & $0.63$   \\
      1      & $0.17$      & $0.27$   & $0.75$   & $0.78$   \\
      2      & $0.20$      & $0.51$   & $0.82$   & $0.83$   \\ \hline
    \end{tabular}
  \end{minipage}
  \begin{minipage}[t]{0.49\textwidth}\centering
    \phantomsubcaption%
    \label{fig:basehdofcol}%
    \begin{tabular}{|c||c|c|c|c|}
      \multicolumn{5}{c}{\subcap{fig:basehdofcol} column-preserving MIS} \\ \hline
      $\ell$ & $\beta=1^0$ & $1^{-2}$ & $1^{-4}$ & $1^{-8}$ \\ \hline
      0      & $0.14$      & $0.14$   & $0.30$   & $0.35$   \\
      1      & $0.17$      & $0.20$   & $0.39$   & $0.47$   \\
      2      & $0.19$      & $0.25$   & $0.48$   & $0.54$   \\ \hline
    \end{tabular}
  \end{minipage}

  \caption{%
    The convergence factor for SA-AMG applied to the \firstblock test problem
    for increasing mesh refinement and decreasing boundary condition strength.
    The aspect ratio of the elements is $\phi=100$; the polynomial
    elements order is $\poly=3$; the meshes are defined by $\ell$ levels
    of isotropic refinement of the mesh $\cT_{xy}$.  We show the convergence
    factor (estimated over the number of iterations needed to reduce the
    residual by a factor of $10^{-14}$ or 100 iterations) for decreasing
    magnitudes of the basal Robin coefficient $\beta$.  Table
    \subcap{fig:basehagg} shows the results for a hierarchy constructed using
    standard randomized MIS aggregate construction, and table
    \subcap{fig:basehdofcol} shows the same for a hierarchy constructed using
    the column-preserving MIS aggregate construction described in the text.%
  }
  \label{fig:baseh}
\end{figure}

\subsection{Construction of low-order approximation $\tilde{F}$}
\label{sec:quadrature}


AMG requires assembled matrices for the construction of
the coarse grid hierarchy. For matrices arising from high-order
discretizations, this assembly requires significant memory and computation
compared to low-order
discretizations \cite{Brown10, HeysManteuffelMcCormickEtAl05}. In three
dimensions, the cost to construct an element matrix $F_i$ is
$\cO(\poly^7)$ per element, or in terms of the number of nodes per element
$\cO(n_{\poly}^{7/3})$.  Instead of the true element matrix $F_i$ for element
$K_i$, we construct an approximation $\tilde{F}_i$ that treats the variables
associated with nodes of the $\tspace_{\poly}(\refel)$ finite element as
variables for a $\poly\times\poly\times\poly$ grid of $\tspace_1(\refel)$
finite elements with corners at the high-order node locations.  This matrix
$\tilde{F}_i$ can be constructed in $\cO(n_{\poly})$ steps.

Aspects of this lower-order preconditioning technique have been studied for
simple problems: spectral equivalence between high-order and lower-order
discretizations of the Laplacian has been proven in two dimensions \cite{Kim07} and is
demonstrated numerically in three dimensions.  The \firstblock $F$ differs from the
operators used in
previous studies in that it involves variable coefficients, high-aspect ratio
elements, and  nonconforming interfaces between elements.  As we will show
below, the interaction of these factors affects the stability and
effectiveness of low-order preconditioning.

\subsubsection{Influence of quadrature on low-order preconditioning}

The effectiveness of the AMG preconditioner for $F$ also
depends on the choice of the quadrature used to construct
$\tilde{F}$.  Let us denote by $\tilde{F}_G$ and $\tilde{F}_{GL}$ the
low-order matrices computed with Gauss-Legendre and Gauss-Lobatto quadrature,
respectively.  We find that $\tilde{F}_{GL}$ leads to a better and more robust
preconditioner, particularly when combined with an incomplete factorization
smoother.
%
The greater stability of Gauss-Lobatto quadrature comes from its
diagonal-lumping behavior \cite{DurufleGrobJoly09}.  In tensor-product
elements, the directional derivative operators that are used to construct the
element matrices can be decomposed as Kronecker products of the form $B \times
B \times D$, where $B$ maps 1D nodal values to 1D values at quadrature points
and $D$ maps 1D nodal values to derivatives at quadrature points.  When we
refer to Gauss-Lobatto quadrature as ``diagonal-lumping'', we mean that for
Gauss-Legendre quadrature, $B$ has off-diagonal entries, whereas for
Gauss-Lobatto quadrature $B$ is the identity.  This does not make
$\tilde{F}_{GL}$ diagonal, but does increase sparsity: $\tilde{F}_{GL}$
contains \abt 30\% fewer nonzero entries than $\tilde{F}_G$.
The factors of an
incomplete factorization can become ill-conditioned when the original matrix
is far from being diagonally dominant, particularly if there is no pivoting
performed during the factorization \cite[Chapter 10.5]{Saad03}.  The
diagonal-lumping inherent in Gauss-Lobatto quadrature increases the
magnitude of diagonal entries relative to off-diagonal entries.
An additional advantage of the sparsity of $\tilde{F}_{GL}$ relative to
$\tilde{F}_G$ is the reduced cost of the incomplete factorization and the
hierarchy construction.

\subsubsection{Influence of nonconforming meshes on low-order preconditioning}
\label{sec:xyz}

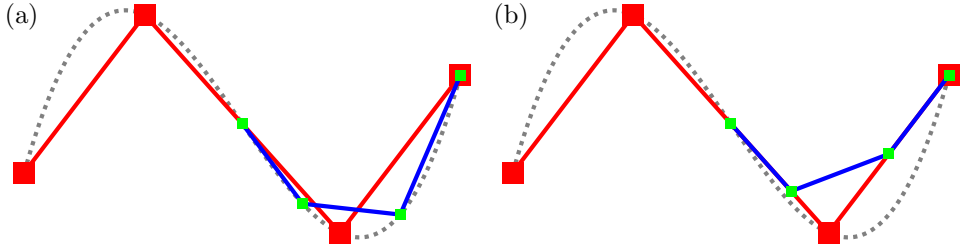
\begin{figure}\centering
  \begin{minipage}[t]{0.49\textwidth}\centering
    \phantomsubcaption%
    \label{fig:interpfull}%
    \makebox[0pt][l]{\subcap{fig:interpfull}}%
    \vtop{%
      \vspace{-0.6\baselineskip}%
      \hbox{%
        \tikzset{interpfull/.style={xscale=2.9,yscale=1.45},
                 interpfull/line/.style=ultra thick}%
        \tikzset{
  interpfull/.prefix style={},
  interpfull/dot/.prefix style={draw=none,rectangle,inner sep=0pt,minimum size=4pt},
  interpfull/bigdot/.prefix style={draw=none,rectangle,inner sep=0pt,minimum size=8pt},
  interpfull/line/.prefix style={},
}
\begin{tikzpicture}[interpfull,
                    dot/.style=interpfull/dot,
                    bigdot/.style=interpfull/bigdot,
                    line/.style=interpfull/line]


  \draw [color=gray,dotted,smooth,samples=200,domain=-1:1,line]
        plot (\x,{\x * (-2.9069 + \x * \x * 3.3541)});

  \draw [draw=red, line]
        (-1.0000cm , -0.4472cm) node[bigdot,fill=red] {} --
        (-0.4472cm ,  1.0000cm) node[bigdot,fill=red] {} --
        ( 0.4472cm , -1.0000cm) node[bigdot,fill=red] {} --
        ( 1.0000cm ,  0.4472cm) node[bigdot,fill=red] {};


  \draw [draw=blue,line]
        ( 0.0000cm  ,  0.0000cm) node[dot,fill=green] {} --
        ( 0.2764cm  , -0.7326cm) node[dot,fill=green] {} --
        ( 0.7236cm  , -0.8326cm) node[dot,fill=green] {} --
        ( 1.0000cm  ,  0.4472cm) node[dot,fill=green] {};

\end{tikzpicture}

%
      }%
    }%
  \end{minipage}
  \begin{minipage}[t]{0.49\textwidth}\centering
    \phantomsubcaption%
    \label{fig:interpsparse}%
    \makebox[0pt][l]{\subcap{fig:interpsparse}}%
    \vtop{%
      \vspace{-0.6\baselineskip}%
      \hbox{%
        \tikzset{interpsparse/.style={xscale=2.9,yscale=1.45},
                 interpsparse/line/.style=ultra thick}%
        \tikzset{
  interpsparse/.prefix style={},
  interpsparse/dot/.prefix style={draw=none,rectangle,inner sep=0pt,minimum size=4pt},
  interpsparse/bigdot/.prefix style={draw=none,rectangle,inner sep=0pt,minimum size=8pt},
  interpsparse/line/.prefix style={},
}
\begin{tikzpicture}[interpsparse,
                    dot/.style=interpsparse/dot,
                    bigdot/.style=interpsparse/bigdot,
                    line/.style=interpsparse/line]


  \draw [color=gray,dotted,smooth,samples=200,domain=-1:1,line]
        plot (\x,{\x * (-2.9069 + \x * \x * 3.3541)});

  \draw [draw=red, line]
        (-1.0000cm , -0.4472cm) node[bigdot,fill=red] {} --
        (-0.4472cm ,  1.0000cm) node[bigdot,fill=red] {} --
        ( 0.4472cm , -1.0000cm) node[bigdot,fill=red] {} --
        ( 1.0000cm ,  0.4472cm) node[bigdot,fill=red] {};


  \draw [draw=blue,line]
        ( 0.0000cm  ,  0.0000cm) node[dot,fill=green] {} --
        ( 0.2764cm  , -0.6180cm) node[dot,fill=green] {} --
        ( 0.7236cm  , -0.2764cm) node[dot,fill=green] {} --
        ( 1.0000cm  ,  0.4472cm) node[dot,fill=green] {};

\end{tikzpicture}

%
      }%
    }%
  \end{minipage}
  \caption{%
    Illustration of different constructions of the low-order matrix $\tilde F$
    at nonconforming interfaces.  $\tilde{F}$ can be thought of as the matrix
    for a finite element space whose functions are piecewise linear between
    the nodes of the original high-order finite element space.  Because these
    nodes do not align at nonconforming interfaces (see \cref{fig:hanging}),
    the functions in this low-order space have discontinuities.  The nature of
    these discontinuities depends on how the nodal values of the smaller
    element, shown in green, are interpolated from the global nodes, which
    coincide with the nodes of the larger element, shown in red.
      \capsubref{fig:interpfull}
        If the same restriction operator $R_i$ is used as for the original
        high-order finite element space, dependent hanging nodes are
        interpolated as though the independent nodes represent high-order
        polynomials.  Each hanging node depends on every independent node
        along the nonconforming edge/face.
      \capsubref{fig:interpsparse}
        If the approximate restriction operator $\tilde{R}_i$ is used,
        dependent hanging nodes are interpolated linearly (or bilinearly for
        two-dimensional faces) from the closest independent nodes.  This
        results in sparser matrices and faster global matrix construction.%
  }%
  \label{fig:interp}%
\end{figure}

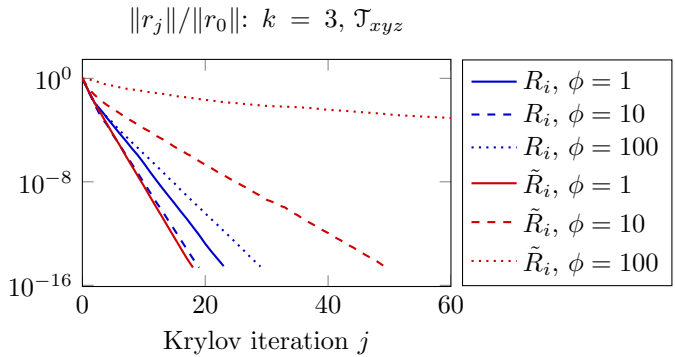
\begin{figure}\centering
  \begin{tikzpicture}
    \begin{axis}
      [
        title      = {$\|r_j\|/\|r_0\|$: $\poly=3$, $\cT_{xyz}$},
        legend pos = outer north east,
        cycle multi list={
          {blue!80!black,thick},{red!80!black,thick}\nextlist
          solid,dashed,dotted
        },
        mark repeat={4},
      ]

      \foreach \j/\name in {full/$R_i$, sparse/$\tilde{R}_i$} 
        {%
          \foreach \i in {1,10,100}
            {%
              \addplot table
                [
                  x index = 0, 
                  y index = 11,
                ]
                {figures/xyz/xyz.\j.\i.monitor};
                \addlegendentryexpanded{\name, $\phi=\i$}
            }
        }
    \end{axis}
  \end{tikzpicture}

  \caption{%
    Relative residual $\|r_j\|/\|r_0\|$ versus number of Krylov iterations for
    the test problems described in \cref{sec:testprob} for the mesh $\cT_{xyz}$
    (\cref{fig:stressxyz}), which has nonconforming faces in all directions.%
    We compare two
    possibilities of handling nonconforming interfaces in the construction of
    the low-order preconditioner $\tilde{F}$ at different element aspect
    ratios.  The first possibility is to use $\{R_i\}$, the same hanging node
    restriction matrices as in the high-order discretization.  An alternative
    is to use $\{\tilde{R}_i\}$, a piecewise linear approximation of
    $\{R_i\}$.%
  }%
  \label{fig:stressconvxyz}%
\end{figure}

As discussed in \cref{sec:fedisc}, $C^0$-conforming discretizations on
nonconforming meshes require element restriction matrices $R_i$ that contain
dense blocks for nonconforming interfaces.  For an element with hanging nodes,
the product $R_i^\tr \tilde{F}_i R_i$ can at best be computed in
$\cO(\poly^5)$ time, or in terms of $n_\poly$, $\cO(n_{\poly}^{5/3})$.  One
possibility for recovering $\cO(n_{\poly})$ element assembly is to replace
$R_i$ with a matrix $\tilde{R}_i$, in which the values of hanging nodes only
depend on the nearest independent nodes, as illustrated in \cref{fig:interp}.
The convergence plots in \cref{fig:base,fig:baseh} were generated on the mesh
$\cT_{xy}$, which is the default mesh for our test problems, using these
sparse-approximation restriction matrices $\{\tilde{R}_i\}$: in this case,
their use does not affect the convergence.   If we use the mesh $\cT_{xyz}$
(\cref{fig:stressxyz}), which is generated by standard isotropic octree-based
refinement, the errors incurred by these sparse restriction operators increase
as $\phi$ increases.  Using the true $\{R_i\}$ for $\cT_{xyz}$ improves the
convergence, but we still observe slight $\phi$-dependence (see
\cref{fig:stressconvxyz}).

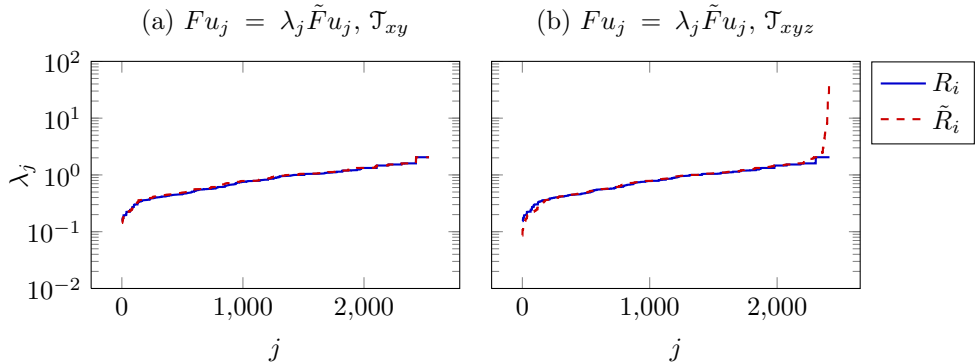
\begin{figure}\centering
  \begin{minipage}{0.47\textwidth}\centering
    \phantomsubcaption%
    \label{fig:diffeigsxy}%
    \begin{tikzpicture}[baseline]
      \begin{axis}
        [
          ytick       = {1.e-2,1.e-1,1.e0,1.e1,1.e2},
          xmin        =,
          xmax        =,
          ymin        = 1.e-2,
          ymax        = 1.e2,
          xlabel      = {$j$},
          title        =%
          {\subcap{fig:diffeigsxy} $Fu_j = \lambda_j \tilde{F}u_j$, $\cT_{xy}$},
          ylabel       = {$\lambda_j$},
          ylabel shift = -1em,
        ]

        \addplot table
          [
            x expr  = \coordindex,
            y index = 0,
          ]
          {figures/diffeigs.xy};

        \addplot table
          [
            x expr  = \coordindex,
            y index = 1,
          ]
          {figures/diffeigs.xy};

      \end{axis}
    \end{tikzpicture}
  \end{minipage}
  \begin{minipage}{0.52\textwidth}\centering
    \phantomsubcaption%
    \label{fig:diffeigsxyz}%
    \begin{tikzpicture}[baseline]
      \begin{axis}
        [
          ytick       = {1.e-2,1.e-1,1.e0,1.e1,1.e2},
          xmin        =,
          xmax        =,
          ymin        = 1.e-2,
          ymax        = 1.e2,
          xlabel      = {$j$},
          yticklabels={},
          title =%
          {\subcap{fig:diffeigsxyz} $Fu_j=\lambda_j \tilde{F}u_j$, $\cT_{xyz}$},
          legend pos = outer north east,
        ]

        \addplot table
          [
            x expr  = \coordindex,
            y index = 0,
          ]
          {figures/diffeigs.xyz};
        \addlegendentry{$R_i$}

        \addplot table
          [
            x expr  = \coordindex,
            y index = 1,
          ]
          {figures/diffeigs.xyz};
        \addlegendentry{$\tilde{R}_i$}
      \end{axis}
    \end{tikzpicture}
  \end{minipage}
  \caption{%
    Comparison of the generalized eigenvalues $\lambda$ that satisfy
    $F\uu=\lambda\tilde{F}\uu$ for two meshes, $\cT_{xy}$ and $\cT_{xyz}$, and for
    $\tilde{F}$ constructed with either the true restriction operator $R_i$ or
    its sparse approximation $\tilde{R}_i$.
      \capsubref{fig:diffeigsxy}
        The generalized eigenvalues for the mesh $\cT_{xy}$ are almost the
        same whether $R_i$ or $\tilde{R}_i$ is used.
      \capsubref{fig:diffeigsxyz}
        The generalized eigenvalues for the mesh $\cT_{xyz}$ are sensitive to
        the use of $\tilde{R}_i$.%
  }%
  \label{fig:diffeigs}%
\end{figure}

To understand this behavior, note that in essence we are comparing two
approximations to the high-order element matrix $R_i^T F_i R_i$: one, $R_i^T
\tilde{F}_i R_i$, where the element matrix is replaced by a low-order
approximation, and another, $\tilde{R}_i^T \tilde{F}_i \tilde{R}_i$, where
additionally the high-order interpolation of the nodal values at nonconforming
faces is replaced by a low-order interpolation. The latter approximation only
affects nonconforming faces; the fact that nonconforming interfaces are much
larger in $\cT_{xyz}$ than in $\cT_{xy}$ for large $\phi$ seems to explain the
different convergence behavior for problems discretized on the two meshes.

To investigate the influence of nonconforming faces on the low-order
preconditioner numerically, we consider the generalized eigenvalue equation
$F\uu=\lambda \tilde{F}\uu$ for $\phi=100$ on the meshes $\cT_{xy}$ and
$\cT_{xyz}$; the results are shown in \cref{fig:diffeigs}.  For eigenvectors $\uu$ with
eigenvalues far from  1, $\tilde{F}$ is a poor approximation of $F$.  For the
mesh $\cT_{xy}$, the errors incurred by the sparse restriction matrices
$\{\tilde{R}_i\}$ due to the nonconforming interfaces are small.  For
$\cT_{xyz}$, errors incurred by $\{\tilde{R}_i\}$ are more significant, and
inspection of the extremal eigenvalues shows that they are associated with
vectors $\uu$ for which $F\uu$ is large only on nonconforming interfaces that
are normal to the $z$-axis.

\subsubsection{Residual computation for smoother on the finest mesh}
\label{sec:residual_smoothing}

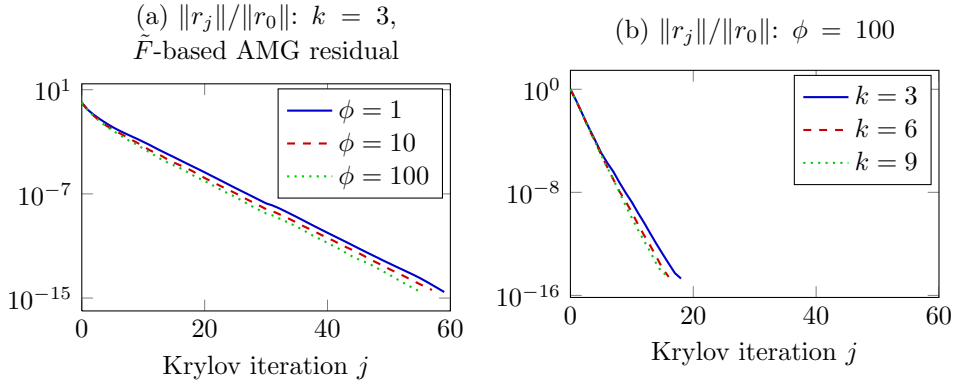
\begin{figure}\centering
  \begin{minipage}[c]{0.49\textwidth}\centering
    \phantomsubcaption%
    \label{fig:pmat}%
    \begin{tikzpicture}[baseline]
      \begin{axis}
        [
          title =%
          {\subcap{fig:pmat} $\|r_j\|/\|r_0\|$: $\poly=3$,\\ $\tilde{F}$-based AMG residual},
        ]

      \foreach \i in {1,10,100} {%
        \addplot table
          [
            x index = 0, 
            y index = 11,
          ]
          {figures/pmat/pmat.\i.monitor};
          \addlegendentryexpanded{$\phi=\i$}
      }
      \end{axis}
    \end{tikzpicture}
  \end{minipage}
  \begin{minipage}[c]{0.49\textwidth}\centering
    \phantomsubcaption%
    \label{fig:poly}%
    \begin{tikzpicture}[baseline]
      \begin{axis}
        [
          title =%
          {\subcap{fig:poly} $\|r_j\|/\|r_0\|$: $\phi = 100$},
        ]

        \foreach \i in {3,6,9} {%
          \addplot table
            [
              x index = 0, 
              y index = 11,
            ]
            {figures/poly/poly.\i.monitor};
            \addlegendentryexpanded{$\poly=\i$}
        }
      \end{axis}
    \end{tikzpicture}
  \end{minipage}
  \caption{%
    Relative residual $\|r_j\|/\|r_0\|$ versus number of Krylov iterations for
    the test problems described in \cref{sec:testprob} on mesh $\cT_{xy}$
    using SA-AMG smoothed by block-Jacobi/IC(0), \subcap{fig:pmat} when the
    residual in the smoother on the finest level of the AMG hierarchy is
    computed with $\tilde{F}$ instead of $F$, and \subcap{fig:poly} for
    varying orders of approximation.
  }%
  \label{fig:stressconv}%
\end{figure}

When the matrix $\tilde{F}$ used to generate a multigrid hierarchy differs
from the true matrix $F$, one has two possibilities for the residual
computation in the smoothing step on the finest mesh. Namely, one can use the
high-order discretized operator, i.e., $r = b - Fx$, or its low-order
approximation, i.e., $r = b - \tilde{F}x$. The convergence in \cref{fig:pmat},
which used $\tilde{F}$ to define the AMG residuals, should be compared to
\cref{fig:stressbase}, which uses $F$ to define residuals and converges
faster.  Thus, in the following, we use the high-order operator for the
residual computation in the smoother on the finest mesh.


\subsection{Convergence for different orders \poly}\label{sec:stresspoly}

In \cref{fig:poly} we show the convergence of our solver (GMRES,
preconditioned by column-preserving SA-AMG with a damped block-Jacobi/IC(0)
smoother) for polynomial orders $\poly$=3, 6, and 9 on a mesh with $\phi=100$.
As can be seen, the iteration number is independent of the polynomial order.
To illustrate that the preconditioner is also $\poly$-independent with respect
to computational work, note that the operator complexities of the AMG
hierarchies with respect to $\tilde{F}$ (i.e., the sum of nonzero matrix
entries in all operators in the AMG hierarchy divided by the number of nonzero
entries in $\tilde F$)  reported by PETSc are 1.22, 1.26, and 1.27 for
$\poly=3,6$ and $9$, respectively, and that the average numbers of nonzeros
per row in $\tilde{F}$ are 86.83, 82.33, and 82.36.  This demonstrates that
the cost to construct $\tilde{F}$ and the coarse hierarchy are proportional to
the problem size, but independent of $\poly$.  The same can be said of the
cost of computing the incomplete factorization smoother on each level as we do
not allow fill-in. Thus, the overall computational complexity for
preconditioning $F$ is independent of $\poly$.

\section{Preconditioning the Schur complement of the \firstblock}
\label{sec:schur}

In this section, we analyze how different approximations of the Schur
complement $S$ of the \firstblock and different choices of the basis for the
pressure space $\cM$ affect the convergence rate for Stokes linear problems
preconditioned as in \cref{eq:ps}.  For all our tests, the preconditioner for
the \firstblock $F$ is a single multigrid V-cycle with the parameter choices
described above in \cref{sec:fblock}.  We test the effectiveness of our
preconditioner on Stokes problems whose setup is discussed in
\cref{sec:testprob}.
As discussed in \cref{sec:fedisc}, the discrete inf-sup constant for the
$\tspace_{\poly} \times \pspace_{\poly - 1}$ mixed element is
$\phi$-dependent.  In \cref{fig:checkplot}, we demonstrate this
$\phi$-dependence numerically.  Because of this instability we prefer the
$\tspace_{\poly}\times \tspace_{\poly-2}$ mixed element for anisotropic
problems.

For constant viscosity $\visc>0$, $S$ is known to be spectrally equivalent to
the scaled pressure mass matrix $-\mu^{-1}M$.  Because of this equivalence, a
common choice is to approximate $S$ by $-\tilde{M}(\visc^{-1})$, which is a
diagonally-lumped approximation to the $\visc^{-1}$-weighted mass matrix.
This Schur complement approximation does not take into account
the anisotropic part of the 4th-order tensor
$\dvisc(\uu)$ defined in \eqref{eq:rank4}.
As an alternative to the weighted mass matrix, the least-squares
commutator, also known as the {BFBt} preconditioner, has proven to be
a good Schur complement approximation, in particular in the presence of
strongly varying coefficients  \cite{ElmanHowleShadidEtAl06,MayMoresi08}.

We have implemented both of these Schur complement preconditioners, and have
found $-\tilde{M}(|\visc|^{-1})$ to be the most efficient for the problems
targeted in this paper.
One reason for this good performance is that we choose
a basis for the discontinuous finite element pressure space that nearly
diagonalizes the mass matrix $M$ for the pressure space. As a consequence, the
effect of replacing $-M(|\visc|^{-1})$ with its mass-lumped counterpart
$-\tilde{M}(|\visc|^{-1})$ is minimal. Such a basis for $\tspace_{\poly-2}(\refel)$
is given by a Lagrange basis for the
tensor-product Gauss nodes.  For a mapped element $K_i$, the mass matrix
remains nearly diagonal provided the mapping is moderately nonaffine. 
To illustrate the effect of the choice of the basis, in
\cref{fig:badbasis} we compare the convergence when using this Lagrange basis
for the Gauss nodes to the convergence when using the Lagrange basis for
Gauss-Lobatto nodes, which is more commonly used as a basis for
tensor-polynomial finite elements.

\begin{figure}\centering
  \begin{minipage}[t]{0.49\textwidth}\centering
    \phantomsubcaption%
    \label{fig:checkplot}%
    \begin{tikzpicture}[baseline]
      \begin{axis}
        [
          title =%
          {\subcap{fig:checkplot} $\|r_j\|/\|r_0\|$: $\poly = 3$},
          legend columns=2,
          legend style={at={(0.45,-0.4)},anchor=north},
          cycle multi list={
            {blue!80!black,thick},{red!80!black,thick}\nextlist
            solid,dashed
          },
          mark repeat={4},
        ]

        \foreach \j/\name in {Q/$\tspace_{\poly-2}$,P/$\pspace_{\poly-1}$} {%
          \foreach \i in {10,100} {%
            \addplot table
              [
                x index = 0, 
                y index = 11,
              ]
              {figures/checkerboard/checkerboard.\j.\i.monitor};
              \addlegendentryexpanded{\name, $\phi=\i$\:}
          }
        }

      \end{axis}
    \end{tikzpicture}
  \end{minipage}
  \begin{minipage}[t]{0.49\textwidth}\centering
    \phantomsubcaption%
    \label{fig:badbasis}%
    \begin{tikzpicture}[baseline]
      \begin{axis}
        [
          title =%
          {\subcap{fig:badbasis} $\|r_j\|/\|r_0\|$: $\poly = 3$},
        ]

        \foreach \i/\name in {G/Gauss,GL/Gauss-Lobatto} {%
          \addplot table
            [
              x index = 0, 
              y index = 11,
            ]
            {figures/basis/basis.\i.monitor};
            \addlegendentryexpanded{\name}
        }

      \end{axis}
    \end{tikzpicture}
  \end{minipage}
  \caption{%
    Relative residual versus number of Krylov iterations for a Stokes test
    problem, posed on the same domain as the test problems in
    \cref{sec:testprob}, but with a uniformly refined mesh.  The convergence
    obtained with different discretizations or bases for the pressure space
    are compared.
      \capsubref{fig:checkplot}
        Comparison between the pressure spaces $\tspace_{\poly-2}$ and
        $\pspace_{\poly-1}$.
      \capsubref{fig:badbasis}
        Comparison between Gauss points and Gauss-Lobatto points for nodal
        bases of the pressure space $\tspace_{\poly-2}$.%
  }%
  \label{fig:checkerboard}%
\end{figure}

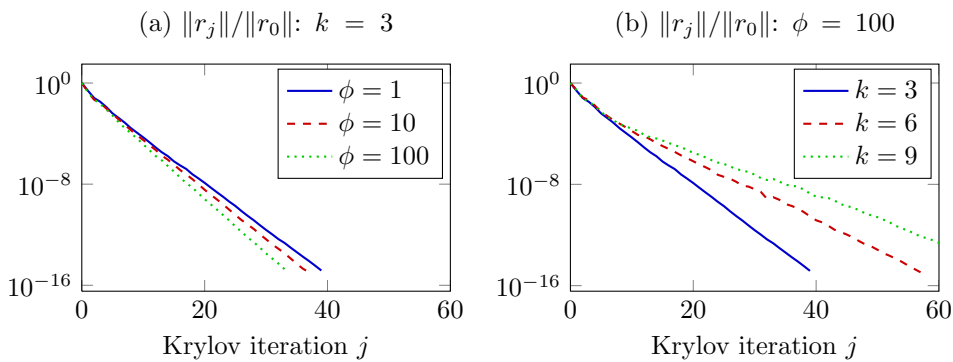
\begin{figure}\centering
  \begin{minipage}[c]{0.49\textwidth}\centering
    \phantomsubcaption%
    \label{fig:basestokes}%
    \begin{tikzpicture}[baseline]
      \begin{axis}
        [
          title =%
          {\subcap{fig:basestokes} $\|r_j\|/\|r_0\|$: $\poly = 3$},
        ]

        \foreach \i in {1,10,100} {%
          \addplot table
            [
              x index = 0, 
              y index = 11,
            ]
            {figures/basestokes/basestokes.\i.monitor};
            \addlegendentryexpanded{$\phi=\i$}
        }

      \end{axis}
    \end{tikzpicture}
  \end{minipage}
  \begin{minipage}[c]{0.49\textwidth}\centering
    \phantomsubcaption%
    \label{fig:polystokes}%
    \begin{tikzpicture}[baseline]
      \begin{axis}
        [
          title =%
          {\subcap{fig:polystokes} $\|r_j\|/\|r_0\|$: $\phi = 100$},
        ]

        \foreach \i in {3,6,9} {%
          \addplot table
            [
              x index = 0, 
              y index = 11,
            ]
            {figures/polystokes/polystokes.\i.monitor};
            \addlegendentryexpanded{$\poly=\i$}
        }

      \end{axis}
    \end{tikzpicture}
  \end{minipage}
  \caption{%
    Relative residual versus number of Krylov iterations for Stokes test
    problems.
      \capsubref{fig:interpfull}
        The convergence for different aspect ratios $\phi$ is compared.
      \capsubref{fig:interpsparse}
        The convergence for different polynomial orders $\poly$ is compared.%
  }%
  \label{fig:stokesplots}%
\end{figure}

In \cref{fig:basestokes}, we observe that the preconditioner for the
Stokes system $A$ results in convergence independent of the aspect ratio
$\phi$; however, as shown in \cref{fig:polystokes}, we find a dependence of the
convergence on \poly. Note that this differs from our findings obtained for
the \firstblock shown in
\cref{fig:poly}, which shows independence of $\poly$ for solving
systems with $F$. The $\mathcal O(\poly^{-1})$ decay of the lower bound for the
discrete inf-sup constant for $\tspace_{\poly} \times \tspace_{\poly-2}$
suggests that this $\poly$-dependence cannot be avoided.

\section{Nonlinear ice stream problems with smooth and rough beds}\label{sec:nonlinex}

\begin{figure}\centering
  \begin{minipage}[t]{0.49\textwidth}\centering
    \phantomsubcaption%
    \label{fig:streambeta}%
    \subcap{fig:streambeta}%
    \vtop{%
      \vspace{-0.6\baselineskip}%
      \hbox{%
        \includegraphics[height=0.68\textwidth]{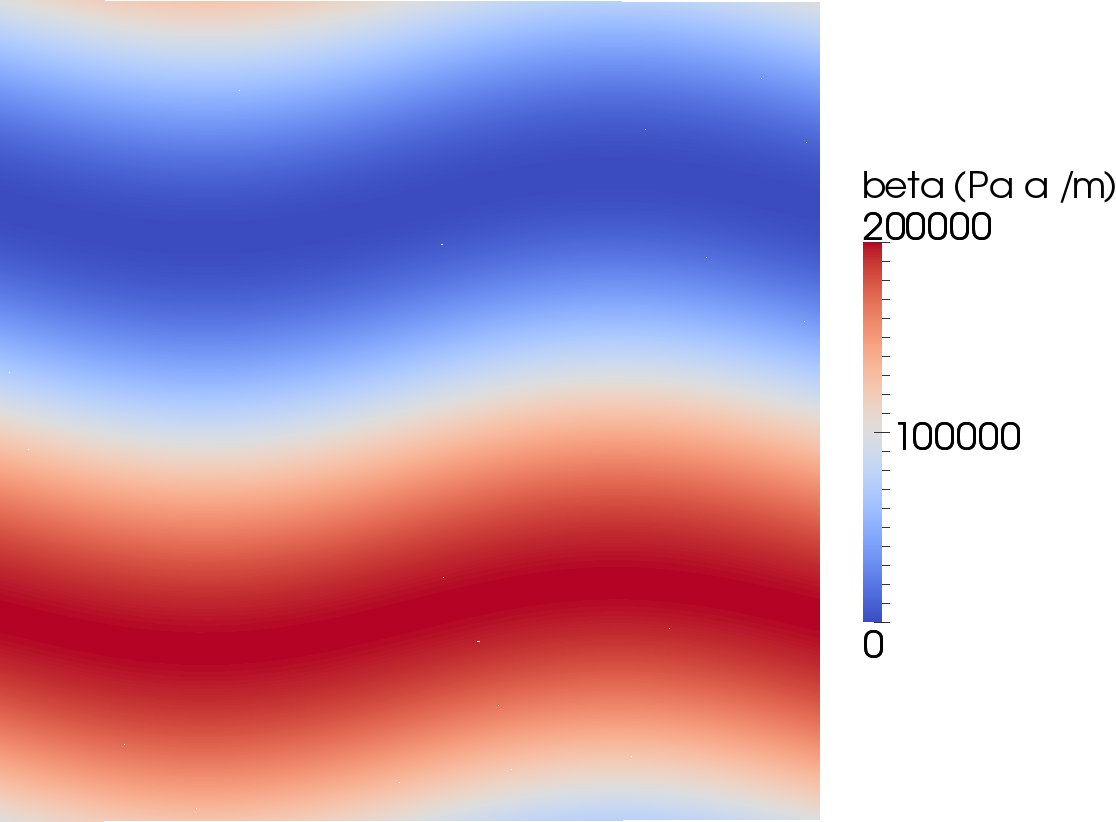}%
      }%
    }%
  \end{minipage}
  \begin{minipage}[t]{0.49\textwidth}\centering
    \phantomsubcaption%
    \label{fig:streamu}%
    \subcap{fig:streamu}%
    \vtop{%
      \vspace{-0.6\baselineskip}%
      \hbox{%
        \includegraphics[height=0.68\textwidth]{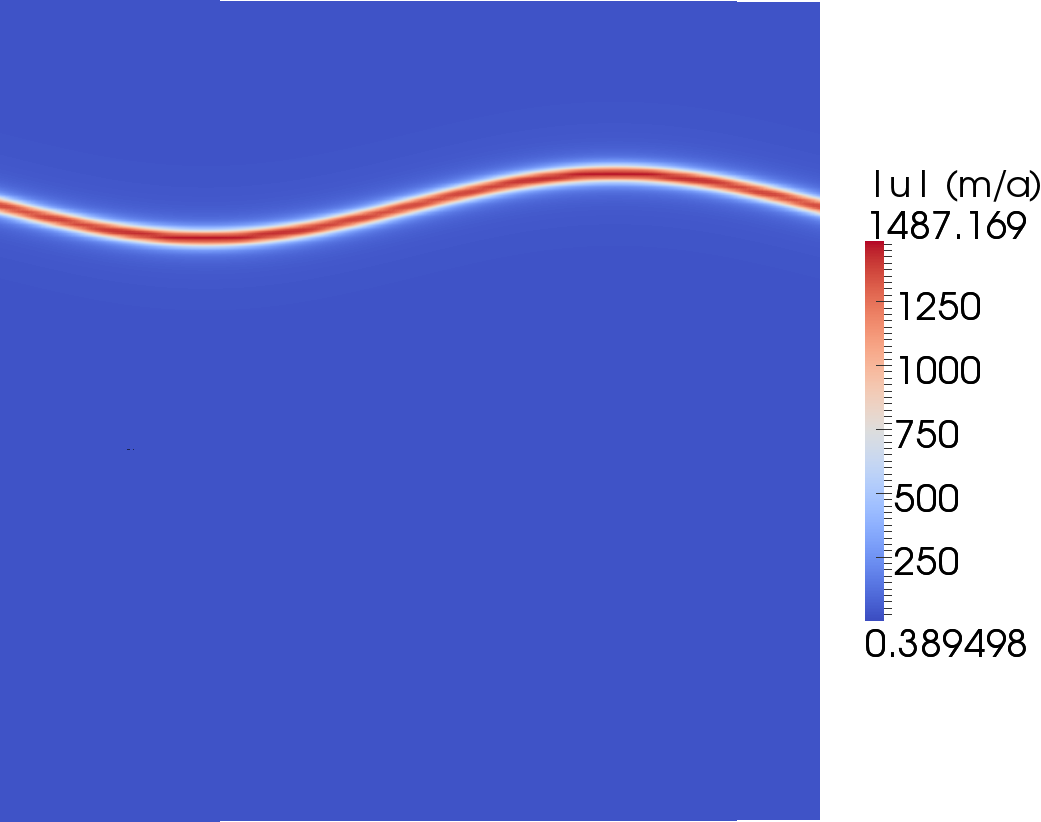}%
      }%
    }%
  \end{minipage}
  \caption{%
    Ice stream model problem on a {\rm{400}$\times$\rm{400} \km} periodic domain of
    {\rm{1} \km} thick ice with a $0.5^\circ$ slope to the right.
      Shown in \capsubref{fig:streambeta} is the
        Robin coefficient field \bta beneath the ice and in
      \capsubref{fig:streamu}
        the magnitude of the velocity \uu at the top surface. A fast flowing ice stream
        develops due to the shear thinning rheology.%
  }%
  \label{fig:stream}%
\end{figure}

To test the nonlinear solver for \cref{eq:weak}, we adapt a model problem from
\cite{CornfordMartinGravesEtAl12}.  As in \cref{sec:testprob}, the domain is a
cutout of an infinite slab that is periodic in the horizontal dimensions, but
the pitch of the domain relative to the direction of gravity is $\theta
=0.5^\circ$, so that a flow is induced.  The Robin coefficient field \bta is
shown in \cref{fig:streambeta}.  Although \bta varies smoothly, the
nonlinearity of the rheology causes the velocity \uu to develop a narrow
region of fast flow similar to an ice stream, as shown in \cref{fig:streamu}.
The constants in the constitutive relationship \cref{eq:glen} are $n=3$ and
$B(T)\equiv 2.15 \times 10^5\ \Pa\ \yr^{1/3}$, which equals $A^{-1/3}$, where
$A=10^{-16}\ \Pa^{-3}\ \yr^{-1}$, which is taken from
\cite{PattynPerichonAschwandenEtAl08}.  Because the top surface of the
periodic domain is flat, we can convert the body force due to gravity into a
constant tangential surface force with magnitude $\rho g \sin\theta$:  we
use an ice density of $\rho = 910 \kg \m^{-3}$, and the acceleration due to
gravity is $g = 9.81 \m \sec^{-2}$.  Note that with this change in the
forcing, the pressure $p$ is now interpreted as the variation from hydrostatic
pressure.  We use $\varepsilon = 1\times
10^{-6}\ \yr^{-2}$, which has a negligible effect for stresses of $10^5\ \Pa$
or greater, following the recommendation in \cite[Chapter 2]{Hutter83}.  The
periodic domain is 400 \km $\times$ 400 \km $\times$ 1 \km, and we again use
$\cT_{xy}$ as our mesh so that the elements are stretched to $\phi=100$.

\begin{figure}\centering
  \begin{tikzpicture}[baseline]
    \begin{axis}
      [
        title = {$\|r_j\|/\|r_0\|$: $\phi = 100$, ISMIP-C, stream variant, \\
        CG(2)/block-Jacobi/IC(0)},
        title style = {align = center, text width=9cm},
        xmax  = 80,
        width = 7.5cm,
        xlabel = cumulative Krylov iteration $j$,
        legend entries={{$\poly=3$, Linear},
                        {$\poly=3$, Newton [9]},
                        {$\poly=4$, Linear},
                        {$\poly=4$, Newton [10]},
                        {$\poly=5$, Linear},
                        {$\poly=5$, Newton [10]},
                        {$\poly=3$, Linear},
                        {$\poly=3$, Picard [50]}},
        legend pos = outer north east,
        cycle multi list={
          {blue!80!black},{red!80!black},%
          {green!80!black},{black}\nextlist
          solid,{only marks,mark=diamond*}
        },
      ]

      \foreach \i in {3,4,5,50} {%
        \addplot table {figures/stream/stream.\i.lin};
        \addplot table {figures/stream/stream.\i.nonlin};
      }

    \end{axis}
  \end{tikzpicture}
  
  \caption{%
    Convergence of the inexact Newton solver for the flat bed topography ice
    stream model problem, for different polynomial orders $\poly$.  Diamonds
    correspond to the nonlinear residual at the start of a Newton step.  The
    solid lines show the linear residuals of the Krylov iterations. The dots
    show where the nonlinear residuals are evaluated and the preconditioners
    are recomputed.  The Krylov method is FGMRES(30).  The smoother for the
    \firstblock AMG preconditioner is two iterations of CG, preconditioned by
    IC(0) on each processor block.   As can be observed, the linear and the
    nonlinear convergence rates coincide close to the solution, as desired for
    an inexact Newton-Krylov method. Far from the solution, where the linear
    and nonlinear residuals are very different, the inexactness in the method
    avoids over-solving of the linear systems, and is thus efficient in terms
    of Krylov iterations \cite{EisenstatWalker96}.  The bracketed numbers in
    the legend indicate the total number of nonlinear iterations.%
  }%
  \label{fig:streamconv}%
\end{figure}
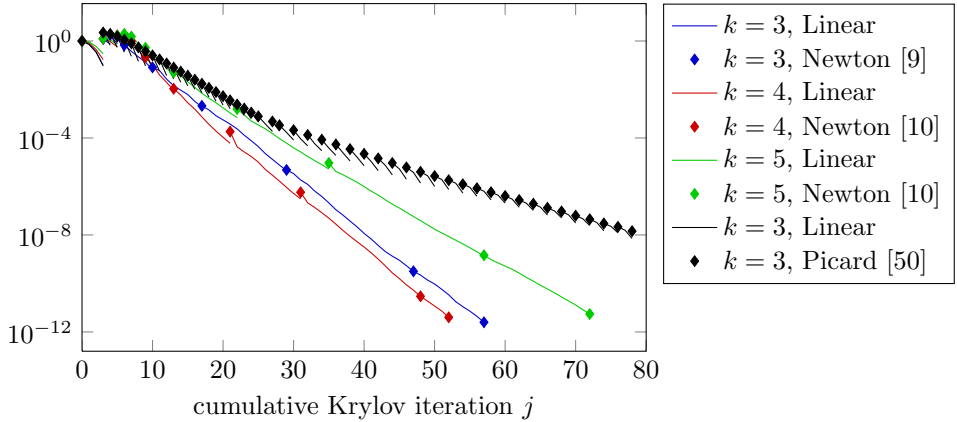

In \cref{fig:streamconv}, we show the convergence behavior of the inexact
Newton method for $\poly=3$, $4$, and $5$, and compare with the convergence of
an inexact Picard method for $\poly=3$.  As can be seen, the inexact Newton
method converges faster than Picard's method, both in terms of the number of
nonlinear iterations and in terms of the total number of Krylov iterations.
As each nonlinear iteration requires a preconditioner setup (or update),  the
superiority of Newton's method compared to the Picard method is even more
pronounced if we consider time-to-solution rather than the number of Krylov
iterations. Although not shown in \cref{fig:streamconv}, we did test
Picard's method with tighter tolerances on the linear solves, but found that
this did not improve the convergence rate in terms of total Krylov iterations.

\begin{figure}\centering
  \begin{minipage}[b]{0.32\textwidth}\centering
    \phantomsubcaption%
    \label{fig:streambeda}%
    \subcap{fig:streambeda}%
    \vtop{%
      \vspace{-0.6\baselineskip}%
      \hbox{%
        \includegraphics[width=0.8\textwidth]%
        {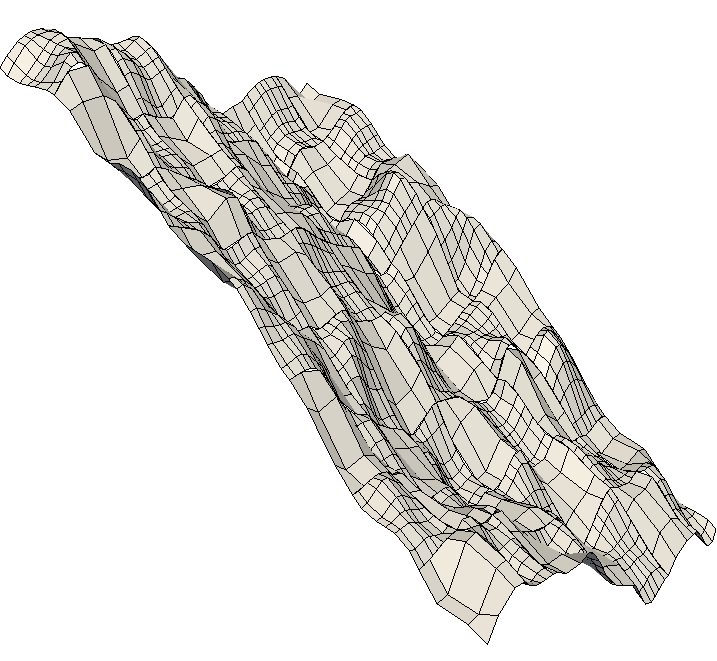}%
      }%
    }%

    \phantomsubcaption%
    \label{fig:streambedc}%
    \subcap{fig:streambedc}%
    \vtop{%
      \vspace{-0.6\baselineskip}%
      \hbox{%
        \includegraphics[width=0.8\textwidth]%
        {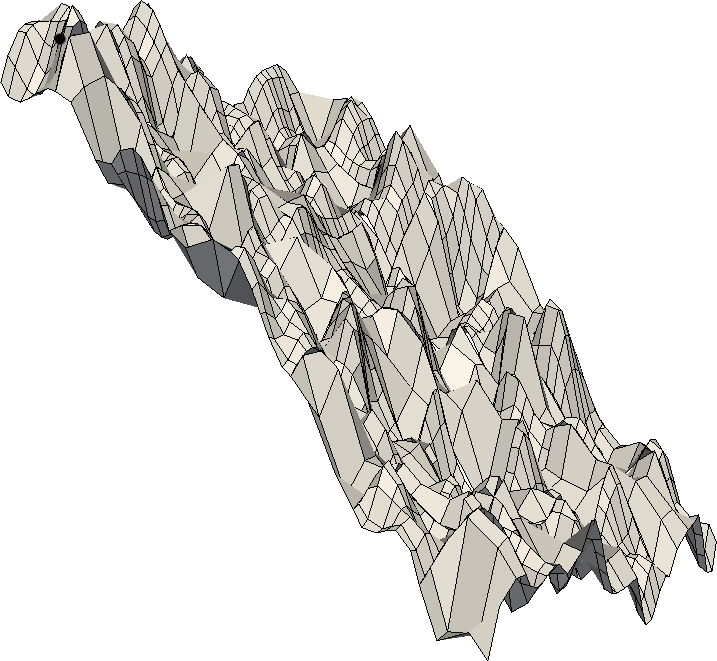}%
      }%
    }%
    \vspace{0.5cm}
  \end{minipage}
  \begin{minipage}[b]{0.65\textwidth}\centering
    \phantomsubcaption%
    \label{fig:bedsconv}%
    \begin{tikzpicture}[baseline]
      \begin{axis}
        [
          title =%
          {\subcap{fig:bedsconv} $\|r_j\|/\|r_0\|$: $\phi = 100$, $\poly = 3$, ISMIP-C, stream
            variant, \\
            CG(2)/block-Jacobi/IC(0)[with shifting ``in blocks''], two bed topographies },
          title style = {align = center, text width=9cm},
          xmax  = 80,
          width = 7.5cm,
          xlabel = cumulative Krylov iteration $j$,
          legend entries={{bed (a)},
                          {},
                          {bed (b)},
                          {}},
          cycle multi list={
            {blue!80!black},{red!80!black},%
            {green!80!black},{black}\nextlist
            solid,{only marks,mark=diamond*}
          },
        ]

        \foreach \i in {1.0,0.5} {%
          \addplot table {figures/stream_wavy/stream.\i.lin};
          \addplot table {figures/stream_wavy/stream.\i.nonlin};
        }

      \end{axis}
    \end{tikzpicture}
    \vspace{-0.1cm}
  \end{minipage}

  \caption{%
    \capsubref{fig:streambeda}, \capsubref{fig:streambedc} Bed topographies
    with different roughness for the ice stream model problem.  The vertical
    scale is exaggerated by a factor of 100. The maximum variation in the bed
    topography in \capsubref{fig:streambedc} is about half of the ice thickness.
    \capsubref{fig:bedsconv} Convergence of the inexact Newton solver for the
    ice stream model problem as described in \cref{fig:stream}, but with 1\%
    of the basal friction $\beta$ and the rough bed topographies.
  }%
  \label{fig:beds}
\end{figure}
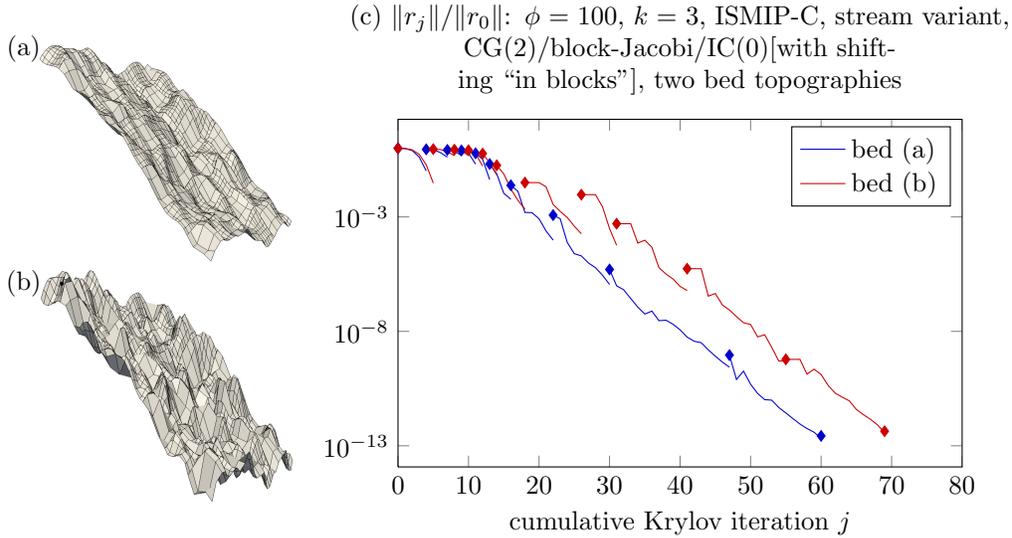

We next test our method on the same problem, but with the Robin coefficient
field reduced to 1\% of the previous field,
and with rough bed geometries instead of the flat slab used in the previous
test. These modifications make the problem more realistic, but more
challenging to solve.  We generate bed topographies using random coefficients
in a truncated Fourier series, as in
\cref{eq:randfield}.  By changing the exponent $\gamma$, which controls the
decay of the Fourier coefficients, we are able to control how rough the
generated topography is.  In \cref{fig:streambeda,fig:streambedc}, we show two topographies,
generated with ten Fourier modes and $\gamma=1.5$ and $\gamma=1.0$,
respectively. In
\cref{fig:bedsconv} we show the convergence behavior of our method on
domains with these bed topographies.

In these nonlinear problems, the incomplete factorization of the \firstblock
approximation $\tilde{F}$ sometimes encounters zero or negative pivots on the
diagonal, which can lead to poor convergence or can cause the solver to fail.
Zero or negative pivots occur more likely when there are regions of rough
topography and sharp solution gradients.  We have considered several strategies
to make the preconditioner robust in these cases.
One remedy to avoid these bad pivots is to increase the sophistication of the
incomplete factorization, for instance by increasing the level of fill-in
(e.g., use IC(1) or IC(2) instead of IC(0)), or by
using a drop-tolerance-based approach. However, these approaches incur
additional
setup and storage costs, and it is difficult to anticipate a~priori
what amount of fill (or what drop tolerance) is needed to avoid bad pivots.

PETSc implements strategies for modification of the matrix being
factorized to avoid bad pivots.  The default is to use Manteuffel's shifting
strategy \cite{Manteuffel80}, in which $\tilde{F}$ is replaced by
$\tilde{F}+\alpha I$.  Our experience is that when the magnitudes of the
entries on the diagonal of $\tilde{F}$ are highly variable, which can be due
to variable coefficients or variable element sizes, this shifting strategy is
detrimental to the effectiveness of the incomplete factorization as a
smoother.  In particular, it can shift the eigenvalues of the high-frequency
modes away from the region that is optimally damped by a Chebyshev polynomial
smoother, leading to stagnation in the convergence.  We find that shifting
``in blocks''---increasing pivots to a lower bound when they become too
small---results more consistently in factorizations that are appropriate for
use as smoothers. The corresponding convergence behavior for problems with
polynomial order $\poly=3$ is shown in \cref{fig:beds}. Note that while
compared to \cref{fig:streamconv}, the efficiency in terms of overall Krylov is
similar, the inexact Newton method requires more iterations before it reaches
the asymptotic convergence regime.  This is likely due to a combination of a
greater degree of nonlinearity in the problem caused by the more complex
geometry and the larger variations in the Robin coefficient $\beta$.

\section{Antarctic ice sheet problem}\label{sec:ant}

We now demonstrate the performance of our inexact Newton-Krylov solver for the
simulation of the dynamics of the Antarctic ice sheet.  Below, we detail how
satellite and radar data are used to define the computational domain and
describe the mesh generation. In \cref{sec:antreggeom}, we study the
performance of our solver for this problem.  Finally, in \cref{subsec:Ant real
geom}, we compare to results presented in \cite{IsaacPetraStadlerEtAl15},
which are obtained with a similar solver that uses a SSOR smoother in the
preconditioner instead of an incomplete factorization smoother.

\subsection{Problem Description}

For the following simulations, we define the ice density to be $\rho =917\ \kg
/ \m^3$, the pre-exponential in Glen's power law \eqref{eq:glen} to be $B(T)=
4.1\times 10^{5}\ \Pa\ \yr^{1/3}$, and the regularizing constant that prevents
infinite effective viscosity to be $\epsilon = 9.95 \times 10^{-6} \yr^{-2}$.
The assumed sliding coefficient field \bta is computed by first taking the
ratio of the driving stress due to gravity and the observed surface velocity,
and then imposing maximum and minimum values.  This results in \bta being
$\abt 0.3\ \MPa\ \yr\ \km^{-1}$ over most of the ice sheet and becomes almost
zero (with a minimum value of $\abt 10^{-14}\ \MPa\ \yr\ \km^{-1}$) in ice
streams.

The body force due to gravity is
$-\rho g\vec z$, where $\vec z$ is a unit vector that is normal to the
reference ellipsoid.
Due to the complexity of the ice sheet geometry, this forcing cannot be
exactly represented by the discrete function spaces we use.
Consequently, we find that when we solve \cref{eq:weak} with the right-hand
side computed by directly discretizing the body force due to gravity, the
velocity in our solution exhibits grid-scale oscillations.  To ameliorate
this, we set the body force to be $\vec b = -\rho g \vec z - \grad \hat{p}$,
where $\hat{p}$ is a hydrostatic-like pressure solving
\begin{align*}
  -\Delta \hat{p} &= 0 &&\vec x \in \Dom,\\
  \hat{p} &= 0 &&\vec x \in \Bndn,\\
  \nabla \hat{p}\cdot \nrml &= -\rho g \vec z\cdot\nrml &&\vec x \in \Bndd,\\
  \nabla \hat{p}\cdot \nrml &= -\rho g \vec z\cdot\nrml &&\vec x \in \Bndr.
\end{align*}
Note that $\hat{p}$ satisfies homogeneous Dirichlet boundary conditions where
the velocity field has stress-free boundary conditions, and
$\hat{p}$ has Neumann boundary conditions where the normal velocity component satisfies a
Dirichlet condition.  The pressure $p$ that is solved for is then a variation
from $\hat{p}$.  Note that the computation of $\hat{p}$ is an upfront
computation and is not included in any of the
performance results that follow.

The geometric description of the ice sheet is constructed from the ALBMAP
dataset \cite{LeBrocqPayneVieli10}.  Elevation values in the ALBMAP dataset
are given relative to the EIGEN-GL04C geoid
\cite{ForsteSchmidtStubenvollEtAl08}: we convert these values to elevations
relative to the WGS84 ellipsoid \cite{USNGA84} using the software library
GeographicLib \cite{Karney13}, and then map the resulting (latitude,
longitude, elevation) geodetic coordinates into Cartesian coordinates.  While
this results in a usable geometry for the ice sheet, the resulting geometry
for the ice shelves, the extension of the ice sheet onto the surface of the
ocean, is far from hydrostatic equilibrium, resulting in flow velocities
several orders of magnitude too fast.  We were unable to correct this behavior
using local geometry adjustments.  While we believe that a good geometry can be
obtained, this was not further pursued for these numerical studies and we have
limited the domain to the grounded ice sheet.

\begin{figure}\centering
  \begin{minipage}[t]{0.49\textwidth}\centering
    \phantomsubcaption%
    \label{fig:antmeshcoarse}%
    \subcap{fig:antmeshcoarse}%
    \vtop{%
      \vspace{-0.6\baselineskip}%
      \hbox{%
        \includegraphics[width=0.9\textwidth]{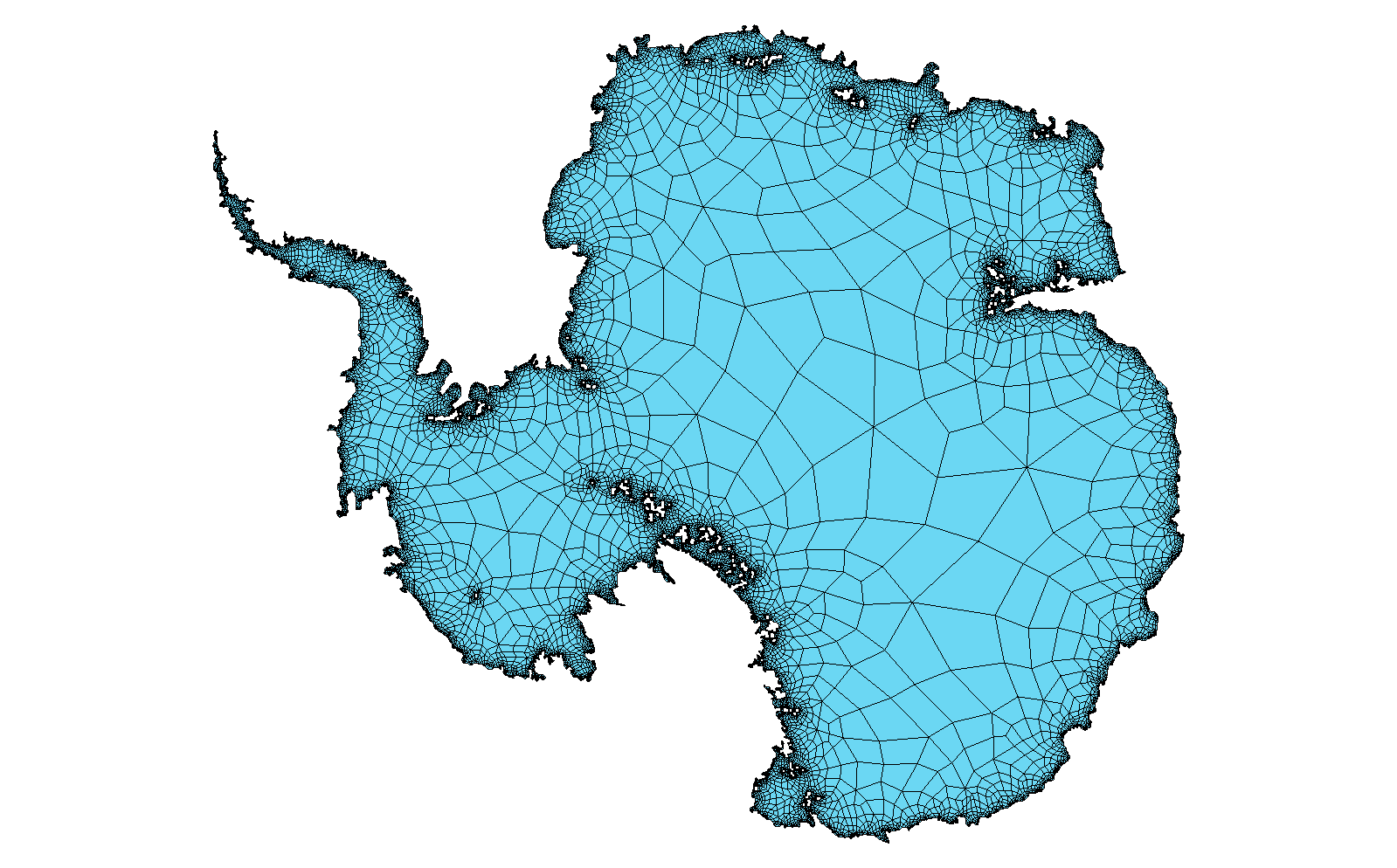}%
      }%
    }%
  \end{minipage}
  \begin{minipage}[t]{0.49\textwidth}\centering
    \phantomsubcaption%
    \label{fig:antmeshdetail}%
    \subcap{fig:antmeshdetail}%
    \vtop{%
      \vspace{-0.6\baselineskip}%
      \hbox{%
        \includegraphics[width=0.9\textwidth]{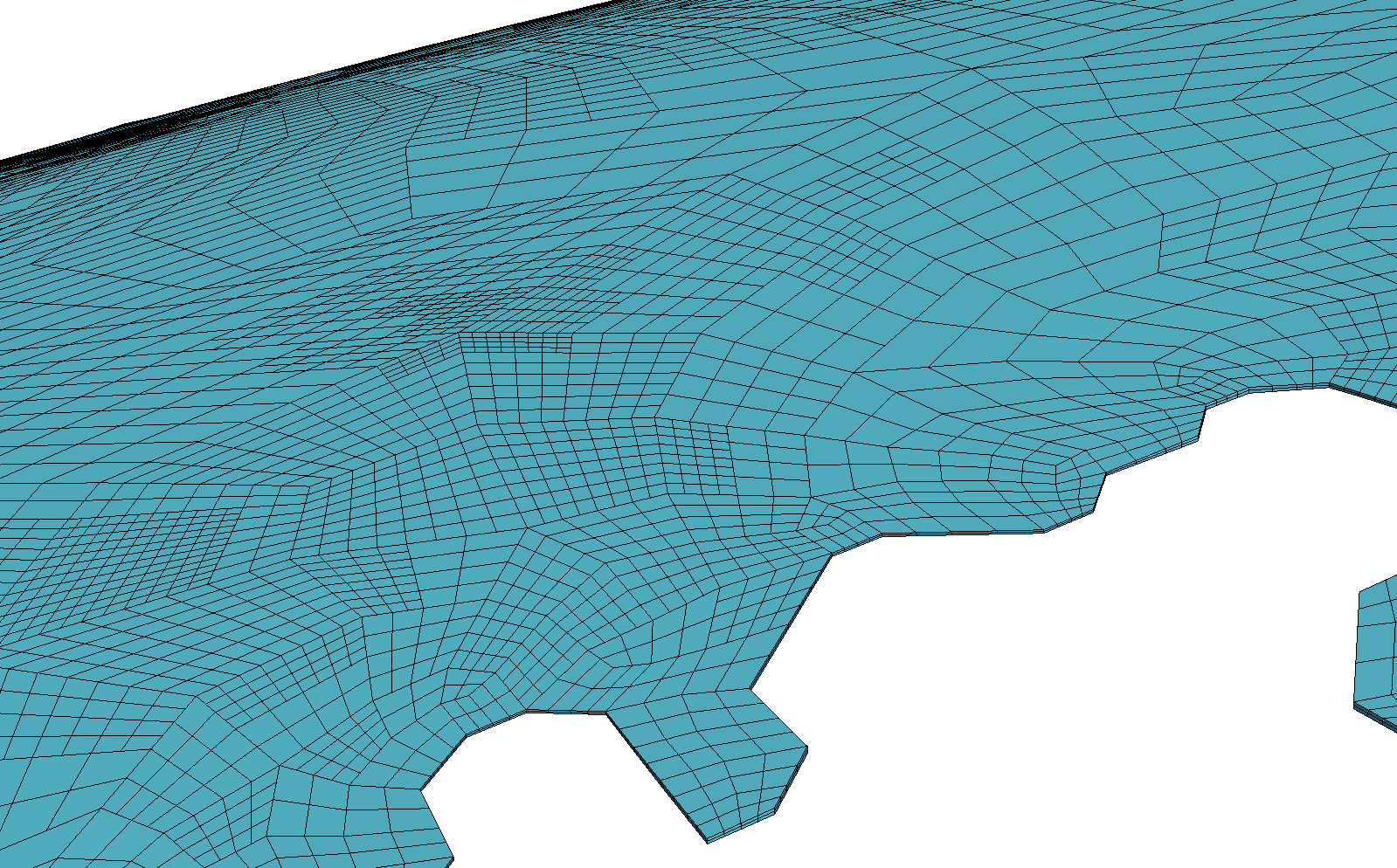}%
      }%
    }%
  \end{minipage}
  \caption{%
    Antarctic ice sheet problem:
      \capsubref{fig:antmeshcoarse}
        The coarse quadrilateral mesh description of the Antarctic ice sheet.
      \capsubref{fig:antmeshdetail}
        An oblique view of the refined hexahedral mesh, obtained using the
        ``\texttt{p6est}'' extension of the \texttt{p4est} library for
        anisotropic AMR.  The refinement near the center is the result of mesh
        refinement to reduce the aspect ratio of elements in the mesh: this
        type of refinement is not possible using purely octree-based
        refinement.%
  }%
  \label{fig:antprob}%
\end{figure}

From the ice thickness data, given on a latitude-longitude grid, we obtain a
polygon describing the lateral boundaries of the ice sheet.  We create a
quadrilateral mesh from this polygon by first using the triangular mesh
generator Triangle \cite{Shewchuk96}, and then splitting the triangles into
quadrilaterals, to which we apply mesh smoothing to improve the element
quality. This resulting coarse quadrilateral mesh (\cref{fig:antmeshcoarse})
contains \abt 27,000 elements.  Using the quadtree-based refinement for the
horizontal directions within our \texttt{p4est} extension for anisotropic
domains, we refine this mesh to construct the footprints for the
columns of our final hexahedral mesh.  We use several refinement criteria: we require that the elements of our mesh have a footprint smaller than
$(2.5\ \km)^2$ at this grounding line; we refine any column whose thickness
varies by more than a factor of 1.5; in keeping with other ice sheet models,
we keep the aspect ratio $\phi$ of the elements below 25.  Once we have
constructed a mesh that satisfies these constraints
(\cref{fig:antmeshdetail}), we use uniform refinement (i.e., replace each
hexahedron by eight children) from this mesh when we perform scaling studies.

The Antarctic ice sheet contains some very thin regions which, because we
constrain the aspect ratio of our elements, would require a large
number of elements if we modeled the true thickness of the ice sheet.  To
control the mesh size, we employ an artificial minimum thickness of 200 \m,
enforced by modifying the bedrock topography.


\subsection{Solver performance and scalability}
\label{sec:antreggeom}
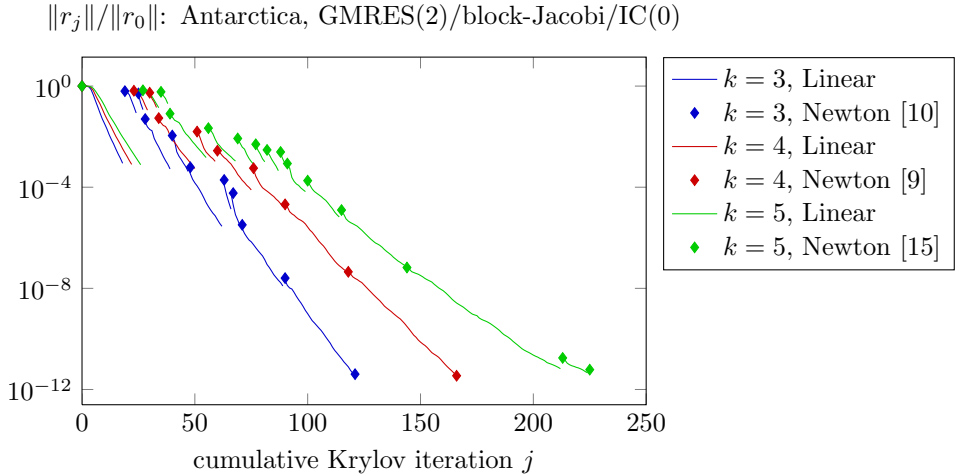
\begin{figure}\centering
  \begin{tikzpicture}[baseline]
    \begin{axis}
      [
        xmax  = 250,
        width = 7.5cm,
        legend entries={{$\poly=3$, Linear},
                        {$\poly=3$, Newton [10]},
                        {$\poly=4$, Linear},
                        {$\poly=4$, Newton [9]},
                        {$\poly=5$, Linear},
                        {$\poly=5$, Newton [15]}},
        legend pos = outer north east,
        xlabel = cumulative Krylov iteration $j$,
        cycle multi list={
          {blue!80!black},{red!80!black},%
          {green!80!black}\nextlist
          solid,{only marks,mark=diamond*}
        },
        title =%
        {$\|r_j\|/\|r_0\|$: Antarctica, GMRES(2)/block-Jacobi/IC(0)},
        title style = {align = center, text width = 9cm},
      ]

      \foreach \i in {3,4,5} {%
        \addplot table {figures/antarctica4901971/ant.\i.stdout.lin};
        \addplot table {figures/antarctica4901971/ant.\i.stdout.nonlin};
      }
    \end{axis}
  \end{tikzpicture}

  \caption{%
    Antarctic ice sheet problem: convergence of the incomplete
    Newton-Krylov solver for different polynomial element orders $\poly$.
    We enforce a minimal ice thickness of 200\ \m, and the mesh is distributed
    over 1024 processes. The bracketed numbers are the total number of
    nonlinear iterations.%
  }%
  \label{fig:antconv}%
\end{figure}

\paragraph{Nonlinear convergence for different orders \poly}
We first test our solver for different polynomial orders on a mesh 
obtained by one level of uniform mesh
refinement. For smoothing in the \firstblock multigrid preconditioner, we use
2 GMRES iterations of our block-Jacobi/IC(0) smoother.  We discretize
\cref{eq:weak} using $\tspace_k \times \tspace_{k-2}$ elements for $\poly=3$,
$4$, and $5$, resulting in problems with 51M, 121M, and 238M degrees of
freedom. The experiments were conducted on TACC's Stampede supercomputer, with
each MPI process assigned to one Sandy Bridge Xeon core.  Each discretization
is distributed across 1024 MPI processes.
The convergence of our inexact Newton-Krylov solver is shown in
\cref{fig:antconv}. Note that the convergence is similar to the results for
the rough bed topography model problem shown in
\cref{fig:bedsconv}. The increased number of overall Krylov iterations (by
about a factor of 2) is likely due to the more complex geometry and boundary
conditions.

\paragraph{Parallel scalability}
\label{sec:parallel}

We next test the parallel scalability of our solver on the $\tspace_3 \times
\tspace_1$ discretization.  We use three different levels of 
uniform mesh refinement, $\ell=0$, $1$, and $2$.  In conducting our study, we found
that our GMRES(2)/block-Jacobi/IC(0) smoother for the \firstblock multigrid
preconditioner was sometimes insufficient: because of the non-smooth
coefficients obtained for some Newton steps, the prolongation of coarse
solutions onto fine grids  would introduce errors that could not be
sufficiently damped by just a couple of smoother iterations.  For this reason,
we employ a stronger GMRES(10)/block-Jacobi/IC(0) smoother in this scaling
study.
We report the results  in \cref{tab:totalstrong}.   The
table includes timings for each of the main components of the nonlinear
solver: the matrix-vector product, the preconditioner application, and the
preconditioner reconstruction.
%

\begin{table}\centering
  \setlength{\parskip}{1.0ex}
  \caption{%
    Antarctic ice sheet problem: scaling of the full nonlinear solver for
    $\poly=3$ for the geometry with minimal ice thickness of 200\ \m, using
    GMRES(10)/block-Jacobi/IC(0) as a smoother in the \firstblock multigrid
    preconditioner.
    For different numbers of processes $P$, we report the total time (in
    seconds) to solve the nonlinear problem to a relative tolerance of
    $10^{-12}$ in the $\ell^2$-norm from a zero initial guess, the parallel
    efficiency (eff.) as well as the number of Newton iterations (\#N) and
    overall Krylov iterations (\#K) performed during the solution.  We also
    report the average time and efficiency for a single Stokes operator
    application (Op), a preconditioner application (PC), and the total
    preconditioner setup time (setup).  Factorization of the coarsest operator
    (which has always less than 50 unknowns) is included in the setup times,
    and applying the inverse of the factored coarse operator is included in
    preconditioner application times. We report efficiency in reference to the
    smallest problem size on the smallest number of processors.
  }
  \begin{tabular}{|r||cc|c|c||cc||cc||cc|} \hline
    $P$   & Solve & eff. & \#N & \#K & Op & eff. & PC & eff. & Setup & eff. \\ \hline
    \multicolumn{11}{|c|}{$\ell=0$ uniform refinement, $\ndof=$7M}               \\ \hline
    128   & 67.7  & 1.00 & 7   & 66  & 0.0275 & 1.00 & 0.981 & 1.00 & 7.361 & 1.00 \\
    256   & 36.9  & 0.91 & 7   & 67  & 0.0146 & 0.94 & 0.577 & 0.93 & 4.255 & 0.86 \\
    512   & 20.3  & 0.83 & 7   & 65  & 0.0080 & 0.86 & 0.299 & 0.82 & 2.722 & 0.67 \\ \hline
    \multicolumn{11}{|c|}{$\ell=1$ uniform refinement, $\ndof=$51M}              \\ \hline
    1,024  & 78.0  & 0.86 & 11  & 75  & 0.0276 & 0.99 & 0.987 & 0.99 & 8.657 & 0.85 \\
    2,048  & 44.2  & 0.76 & 10  & 75  & 0.0152 & 0.90 & 0.561 & 0.87 & 5.810 & 0.63 \\
    4,096  & 30.2  & 0.56 & 10  & 75  & 0.0087 & 0.79 & 0.383 & 0.64 & 5.406 & 0.34 \\ \hline
    \multicolumn{11}{|c|}{$\ell=2$ uniform refinement, $\ndof=$383M}             \\ \hline
    8,192  & 108.0 & 0.62 & 13  & 91  & 0.0295 & 0.93 & 1.13  & 0.87 & 15.58 & 0.47 \\
    16,384 &  74.7 & 0.45 & 10  & 89  & 0.0168 & 0.82 & 0.80  & 0.61 & 17.22 & 0.21 \\ \hline
  \end{tabular}
  \label{tab:totalstrong}%
\end{table}


The time to apply the preconditioner is almost entirely spent in applying the
AMG V-cycle to the \firstblock of the Stokes operator. This step also requires more
communication than the matrix-vector product, both to project and restrict
vectors in the hierarchy and to apply the smoothers.  Note that GMRES
requires global reductions to compute the required inner products, which would
not be the case for a stationary smoother (e.g. Chebyshev).

It is well known that the setup phase in parallel implementations of algebraic
multigrid requires significant communication to properly aggregate degrees of
freedom across processor boundaries, compute prolongation matrices and
repartition coarse matrices.  As a consequence, the setup is often the
least scalable component of the solver \cite{ChowFalgoutHuEtAl06,
  BursteddeGhattasStadlerEtAl09}, which can also be seen in
\cref{tab:totalstrong}.  Here, the reported times include both the prolongator
constructions in the initial setup and the Galerkin projections to reconstruct
the coarse matrices during subsequent Newton steps, so the increase in Newton
iterations for the larger problems affects the reported efficiency.  In all
cases, the contribution of the initial prolongator setup is roughly equivalent
to the cost of one subsequent preconditioner reconstruction.

Note that though the total number of Krylov iterations increases by $\abt37\%$
from the smallest to the largest problem, the number of outer inexact Newton
steps almost doubles from 7 to 13.  Recent work \cite{BrownSmithAhmadia13} on
solving the hydrostatic approximation to the Stokes equations for ice sheet
dynamics has demonstrated the effectiveness of grid continuation in obtaining
an initial guess that is near the region of asymptotic convergence of Newton's
method.  The hierarchical mesh refinement we use lends itself naturally to
grid continuation, so this approach could improve the efficiency of our
nonlinear solver.

\subsection{Comparison with performance of SSOR-based smoothing}

\label{subsec:Ant real geom}

\begin{figure}\centering
          \includegraphics[width=0.9\textwidth]%
          {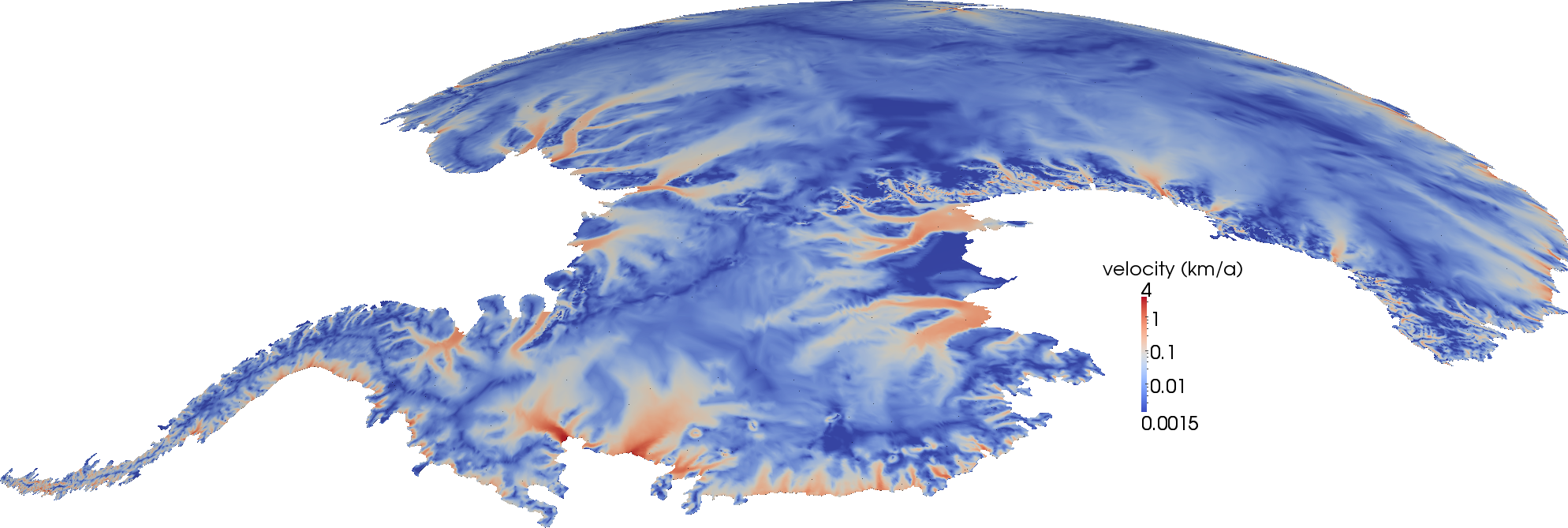}
  \caption{%
    Antarctic ice sheet problem:
    magnitude of the surface velocity field optimized to match satellite observations.
    The inversion procedure to infer the basal coefficient is described in~\cite{IsaacPetraStadlerEtAl15}.
  }%
  \label{fig:invresults}%
\end{figure}

In \cite{IsaacPetraStadlerEtAl15}, we have used a similar nonlinear and linear
solver framework as the one developed in this paper to infer the Robin
coefficient field $\beta$ in the Stokes boundary value problem such that the
velocity fields of the solution closely matches satellite observations of the
Antarctic ice sheet's surface (see \cref{fig:invresults}).  This inference
uses methods of adjoint-based PDE-constrained optimization to find an optimal
$\beta$, which requires the solution of the Stokes problem described in this
work, as well as linearized adjoint problems, whose operators are similar to
the linearized Stokes operator in \eqref{eq:newton}.

The difference between the solver used in \cite{IsaacPetraStadlerEtAl15} and
the one presented here lies in the preconditioning for the \firstblock. In
\cite{IsaacPetraStadlerEtAl15}, a Chebyshev(2)/block-Jacobi/SSOR smoother
combined with the standard SA-AMG scheme (not the column-preserving SA-AMG
presented in \cref{sec:amg}) is used.  In \cref{tab:titan}, we reproduce some
of the algorithmic scalings from \cite{IsaacPetraStadlerEtAl15}, and compare
them to scalings obtained for the same problem with the solvers developed
here: using a Chebyshev(2)/block-Jacobi/IC(0) smoother and column-preserving
SA-AMG. We present only algorithmic scalability because the results in
\cite{IsaacPetraStadlerEtAl15} were computed on Oak Ridge National
Laboratory's Titan supercomputer, while the new scalings were computed on
TACC's Stampede supercomputer. As a consequence, the timings are not directly
comparable.  We note, however, that the cost of applying an IC(0)
preconditioner, in terms of memory movement and floating point operations, is
almost identical to the cost of applying an SSOR preconditioner, so the only
significant additional computation in our new scaling results is in the
computation of the IC(0) factorization, which is on par with a single
application of the IC(0) preconditioner.  Thus, the algorithmic speedup shown
in \cref{tab:titan} is representative of the runtime speedup we would find,
had these scaling results been produced on the same architecture.

We note that the total number of iterations required to solve these problems
is lower than for the problems solved earlier in \cref{sec:parallel}, even
though both are posed on the same domain.  This is because the aspect ratio of
the elements is smaller, with all element aspect ratios $\phi$ less than 10,
and because the basal friction is larger, reaching $\abt3\ \MPa\ \yr\
\km^{-1}$ in some areas (versus the $\abt0.3\ \MPa\ \yr\ \km^{-1}$ used in
\cref{sec:parallel}), and having a minimum value in the ice streams that is three
orders of magnitude larger than the minimum used in \cref{sec:parallel}.

\begin{table}
  \caption{%
    Antarctic ice sheet problem: comparison of algorithmic scalability of the
    solver used in \cite{IsaacPetraStadlerEtAl15} (which uses SSOR and
    standard SA-AMG within the \firstblock preconditioner) versus the one
    presented here (which uses IC(0) and column-preserving SA-AMG).  Shown are
    the number of overall degrees of freedom $\ndof$, the number of processes
    $P$, the number of Newton iterations $\#N$ and the number of cumulative
    Krylov iterations $\#K$.  The data on the left is reproduced from
    \cite{IsaacPetraStadlerEtAl15}.%
  }
  \begin{center}
    \renewcommand{\arraystretch}{1.1}
    \begin{tabular}{|c|r||c|c||c|c|}
      \hline
      \multicolumn{2}{|c||}{} &
      \multicolumn{2}{c||}{SSOR, standard SA-AMG} &
      \multicolumn{2}{c|}{IC(0), column-preserving SA-AMG} \\ \hline
      $\ndof$ & $P$  & \#N & \#K & \#N & \#K \\ \hline
      38M     & 1,024 & 8   & 147 & 7   & 85  \\ \hline
      268M    & 8,192 & 9   & 243 & 8   & 98  \\ \hline
    \end{tabular}
  \end{center}
  \label{tab:titan}

\end{table}

\section{Conclusions}
\label{sec:discussion}

Several issues related to high-order, adaptive mesh discretizations and
solvers for the simulation of nonlinear Stokes flow in three-dimensional
anisotropic domains are addressed in this work. Our main target problem has
been the nonlinear Stokes boundary value problem arising in ice sheet
dynamics. 
We demonstrate an extension to the \texttt{p4est} library for adaptive mesh
refinement of anisotropic domains that combines quadtree-based refinement in
the horizontal directions with columns of elements to achieve a flexible
approach to mesh refinement, with local control over the aspect ratio of
elements in the mesh.  We demonstrate that high-order finite elements
discretized on these meshes are well-approximated by low-order approximations,
which can be used for preconditioning.  We present an efficient solver for the
linearized Stokes equations, with particular emphasis on the design of
algebraic multigrid solvers for high-order discretizations, anisotropic
domains, and hanging nodes.  Using incomplete factorization-based smoothing
for the \firstblock yields efficient and fast convergence.  When using this
type of smoothing, we demonstrate that column-preserving SA-AMG, as
implemented by our DofColumns plugin for PETSc's GAMG preconditioner,
significantly improves over standard SA-AMG in its effectiveness, especially in
the presence of weak boundary conditions. The numerical experiments on our
discretization of the Antarctic ice sheet show that, up to a point, the
incomplete factorization process can be made robust against bad pivots that
occur in problems with variable coefficients and rough topography by using
only local shifting and by combining the incomplete factorization with a
non-stationary Krylov method as a smoother.

\section{Acknowledgments}
\label{sec:acknowledgments}

We would like to thank the three referees who reviewed this work, whose
reviews greatly improved the final version of this work.  In particular, we
would like to thank the two initial referees, Ray Tuminaro and one who remains
anonymous: the column-preserving SA-AMG that we have developed in the
DofColumns plugin was not included in the initial submission of this work, but
was developed in response to their insightful comments and suggestions.

\appendix

\section{Well-posedness of \cref{eq:weak}}
\label{sec:proof}
The main difference between \cref{eq:weak} and the
variational form in \cite{JouvetRappaz12} is the boundary integral $\int_{\Bndr}
\bta\tang\vv\cdot\tang\uu\,d\vec{s}$: for $\uu\in\cW^{1,r}(\Dom)$, this form
is not meaningful because the trace of $\uu$ is not necessarily in
$\cL^2(\Bndr)$.  We compensate for this by bringing the boundary integral
into the definition of $\cV$.
Let $[\cC^\infty(\Dom)]^3_0$ be the space of smooth vector-valued functions in
$\Dom$ that satisfy the homogeneous Dirichlet part of the boundary
conditions and let $r=1 + \frac{1}{n}$.  Under the above assumptions, the
functional
\begin{equation}\label{eq:fnl}
  \cI(\uu) =  \left\{\int_{\Dom} |\Grad\uu|^r\ d\xx\right\}^{1/r} +
  \left\{\int_{\Bndr} \bta |\tang\uu|^2\,d \vec{s}\right\}^{1/2}
\end{equation}
defines a seminorm on $[\cC^\infty(\Dom)]^3_0$. We assume a problem setup in
which $\cI(\uu)$ also defines a norm, which amounts to requiring that if
$\uu\in[\cC^\infty(\Dom)]^3_0$ is a rigid-body motion that satisfies the
Dirichlet conditions, then $\int_{\Bndr}\bta |\tang\uu|^2\,d\vec{s}>0$.  We
define $\cV$ to be the closure of $[\cC^\infty(\Dom)]^3_0$ in this norm.  

Following the same steps as in \cite{JouvetRappaz12}, one can show that the
minimization problem
\begin{equation}
 \min_{\uu}\ \frac{2n}{1+n} \displaystyle\int_{\Dom} B(T) (\IIDu +
 \varepsilon)^{\frac{1 + n}{2n}} d\xx +
  \half \displaystyle\int_{\Bndr} \tang\uu \cdot \tang\uu\ d\vec{s} - \ff(\uu)
\end{equation}
is well-posed in $\cV_{\text{div}}$, the subspace of $\cV$ containing only the
divergence-free functions.  To prove that \cref{eq:weak} is well-posed, we
additionally need to choose a space $\cM$ for which we can prove the inf-sup
condition
\begin{equation}
  \inf_{q\in\cM} \sup_{\uu\in\cV}\ \frac{\int_{\Dom} q \div\uu\
    d\xx}{\|q\|_{\cM}\|\uu\|_{\cV}}\geq \gamma > 0.
\end{equation}
This inequality still holds for $\cM=\cL^{r'}(\Dom)$.  To see this, consider
the subspace $\tilde{\cV}=\{\uu\in\cV:\uu|_{\Bndr} = 0\}$: for $\uu \in
\tilde{\cV}$, $\|\uu\|_{\cV}=\|\uu\|_{[\cW^{1,r}(\Dom)]^3}$, and so
\begin{equation}
  \inf_{q\in\cL^{r'}(\Dom)} \sup_{\uu\in\cV}\ \frac{\int_{\Dom} q \div\uu\
    d\xx}{\|q\|_{\cM}\|\uu\|_{\cV}}\geq \inf_{q\in\cL^{r'}(\Dom)}
  \sup_{\uu\in\tilde{\cV}}\ \frac{\int_{\Dom} q \div\uu\
    d\xx}{\|q\|_{\cL^{r'}(\Dom)}\|\uu\|_{[\cW^{1,r(\Dom)}]^3}}.
\end{equation}
As long as $\Bndn\neq\emptyset$, the term on the right is bounded from below
as a particular case of the inf-sup condition in \cite{JouvetRappaz12}. Thus,
\eqref{eq:weak} is well posed and has a unique solution.

\bibliographystyle{siam}
\bibliography{ccgo}

\end{document}